\documentclass[12pt]{amsart}
\usepackage{amsmath}
\usepackage{amssymb}
\usepackage{amscd}
\usepackage{amsthm}
\usepackage{enumerate}
\usepackage{fullpage}
\usepackage{mathrsfs}
\usepackage{setspace}
\usepackage{mathtools}
\usepackage{xcolor}
\usepackage[all]{xy}
\usepackage{url}
\usepackage[colorlinks=true,linkcolor=blue,breaklinks,pagebackref]{hyperref}
\usepackage[lite,abbrev,alphabetic]{amsrefs}
\usepackage[poorman]{cleveref}
\usepackage{bbm}
\usepackage{tikz-cd}
\usepackage{graphicx}

\theoremstyle{definition}
\newtheorem{definition}{Definition}[section]
\newtheorem{example}[definition]{Example}

\theoremstyle{plain}
\newtheorem{theorem}[definition]{Theorem}
\newtheorem{proposition}[definition]{Proposition}
\newtheorem{lemma}[definition]{Lemma}

\newtheorem*{theorem*}{Theorem}

\theoremstyle{remark}
\newtheorem{remark}[definition]{Remark}

\crefname{theorem}{Theorem}{Theorems}
\crefname{proposition}{Proposition}{Propositions}
\crefname{lemma}{Lemma}{Lemmas}
\crefname{corollary}{Corollary}{Corollaries}
\crefname{conjecture}{Conjecture}{Conjectures}
\crefname{hypothesis}{Hypothesis}{Hypotheses}
\crefname{remark}{Remark}{Remarks}
\crefname{condition}{Condition}{Conditions}
\crefname{example}{Example}{Examples}

\def\C{\mathbb{C}}
\def\R{\mathbb{R}}
\def\Q{\mathbb{Q}}
\def\Z{\mathbb{Z}}
\def\A{\mathbb{A}}

\def\GL{\mathrm{GL}}

\def\inf{\infty}

\def\bs{\backslash}

\DeclareMathOperator{\diag}{diag} 
\DeclareMathOperator{\supp}{supp}
\DeclareMathOperator{\vol}{vol}

\DeclareMathOperator{\rank}{rank}
\DeclareMathOperator{\Hom}{Hom}
\DeclareMathOperator{\Lie}{Lie}

\DeclareMathOperator{\Ad}{Ad}
\DeclareMathOperator{\ad}{ad}

\newcommand{\fa}{\mathfrak{a}}
\newcommand{\fb}{\mathfrak{b}}

\newcommand{\fg}{\mathfrak{g}}

\newcommand{\fm}{\mathfrak{m}}

\newcommand{\fo}{\mathfrak{o}}

\newcommand{\fS}{\mathfrak{S}}

\newcommand{\fI}{\mathfrak{I}}

\newcommand{\bT}{\mathbb{T}}

\newcommand{\cC}{\mathcal{C}}
\newcommand{\cD}{\mathcal{D}}
\newcommand{\cE}{\mathcal{E}}
\newcommand{\cF}{\mathcal{F}}

\newcommand{\cI}{\mathcal{I}}
\newcommand{\cJ}{\mathcal{J}}

\newcommand{\cL}{\mathcal{L}}

\newcommand{\cN}{\mathcal{N}}
\newcommand{\cO}{\mathcal{O}}
\newcommand{\cP}{\mathcal{P}}
\newcommand{\cQ}{\mathcal{Q}}
\newcommand{\cR}{\mathcal{R}}
\newcommand{\cS}{\mathcal{S}}

\newcommand{\cU}{\mathcal{U}}

\numberwithin{equation}{section}

\newcommand{\tG}{\widetilde{G}}
\newcommand{\tM}{\widetilde{M}}
\newcommand{\tL}{\widetilde{L}}
\newcommand{\tP}{\widetilde{P}}
\newcommand{\tQ}{\widetilde{Q}}
\newcommand{\tR}{\widetilde{R}}

\newcommand{\tS}{\widetilde{S}}

\newcommand{\talpha}{\widetilde{\alpha}}
\newcommand{\tbeta}{\widetilde{\beta}}

\newcommand{\tvarpi}{\widetilde{\varpi}}
\newcommand{\tsigma}{\widetilde{\sigma}}

\newcommand{\tfo}{\widetilde{\mathfrak{o}}}
\newcommand{\tcO}{\widetilde{\mathcal{O}}}
\newcommand{\tO}{\widetilde{O}}
\newcommand{\tf}{\tilde{f}}
\newcommand{\sC}{\mathscr{C}}
\newcommand{\sU}{\mathscr{U}}

\newcommand{\sat}{\mathrm{sat}}
\newcommand{\level}{\operatorname{level}}
\newcommand{\fO}{\mathfrak{O}}

\newcommand{\tomega}{\widetilde{\omega}}

\begin{document}

\title{On the continuity of the geometric side of the twisted trace formula}

\author{Satoshi Wakatsuki}
\address{Faculty of Mathematics and Physics, Institute of Science and Engineering, Kanazawa University, Kakumamachi, Kanazawa, Ishikawa, 920-1192, Japan}
\email{wakatsuk@staff.kanazawa-u.ac.jp}

\date{\today}

\begin{abstract}
We establish the continuous extension of the geometric side of the twisted trace formula to a natural space $\cC(\tG(\A)^1)$ of test functions, which are not necessarily compactly supported, for reductive algebraic groups. 
The continuity of the geometric side of the (untwisted) Arthur trace formula was established by Finis--Lapid \cites{FL11semi,FL16}. 
Subsequently, Parab \cite{Par19} generalized their result to the twisted trace formula under a certain restriction on the truncation parameter and an additional assumption.
In this paper, we remove his restriction and assumption, and prove the continuous extension for the twisted trace formula. 
We also give a revised form of Parab's decay estimate for the geometric side.
\end{abstract}

\maketitle

\tableofcontents

\section{Introduction}\label{sec:intro}

Arthur's trace formula is a fundamental tool in the theory of automorphic forms, relating the spectrum of automorphic
representations to the geometry of conjugacy classes.
Its twisted form has played an important role in the study of automorphic representations, notably in base change and endoscopy.
For analytic applications of the trace formula, it is desirable to treat both sides on spaces of test functions that need not have compact support.
This requires not only the continuity of the individual terms in the expansions, but also the absolute convergence of their sums and integrals, together with uniform estimates.
In the untwisted case, following their work on the semisimple terms \cite{FL11semi}, Finis--Lapid established a continuous extension of the entire geometric side to a natural space of functions all of whose derivatives belong to $L^1$; see \cite{FL16}.
Parab \cite{Par19} obtained a twisted analogue for truncation parameters lying on a fixed line.
His argument relies on the root cone lemma, a statement about root systems that he proved for split groups.
In this paper, we establish a continuous extension of the twisted geometric side for general reductive algebraic groups, including nonsplit groups, without assuming this lemma or restricting the truncation parameter to a line.
We also revise Parab's decay estimate for the geometric side to account for additional contributions arising from the twist.

Throughout this paper, the ground field is $\Q$, and $\A$ denotes its ring of adeles.
Let $G$ be the identity component of a reductive algebraic group $G'$ over $\Q$, and fix a connected component $\tG$ having rational points.
For an open compact subgroup $K_0\subset G(\A_f)$, let $\cC(\tG(\A)^1;K_0)$ denote the space of smooth,
right-$K_0$-invariant functions whose right derivatives of every order at the archimedean place belong to $L^1$.
We equip this space with the Fr\'echet topology defined by the $L^1$ norms of these derivatives and write
$\cC(\tG(\A)^1)$ for the inductive limit over $K_0$ (see \S~\ref{sec:Frechet} for details). 
The test functions $f$ in this space $\cC(\tG(\A)^1)$ are not required to have compact support.

The main results are Theorems~\ref{thm:main theorem 1} and~\ref{thm:main theorem 2}.
In the coarse expansion of Theorem~\ref{thm:main theorem 1}, elements are grouped according to the rational conjugacy classes of their semisimple parts. 
Theorem~\ref{thm:main theorem 2} refines each coarse class according to the geometric conjugacy classes of the elements themselves (see \S~\ref{sec:mainrefinedclass} for the definition). 
We state the conclusions for the refined expansion and write $\cO$ for the set of refined classes.
Let $d(T)=\min_{\alpha\in\Delta_{P_0}}\alpha(T)$.
For every $T\in\fa_{M_0}^G$ satisfying $d(T)>d_*$, where $d_*$ is independent of $K_0$ and $f$, the sum over $\fo\in\cO$ of the $L^1$ norms of the truncated kernels is finite.
Consequently, the expansion
\[
 J^{T,\tG}(f)=\sum_{\fo\in\cO}J_{\fo}^{T,\tG}(f)
\]
converges absolutely.
The seminorm bound in Theorem~\ref{thm:main theorem 2} also shows that each $J_{\fo}^{T,\tG}$ is a continuous linear form on $\cC(\tG(\A)^1;K_0)$.
For each $f$ and each $\fo\in\cO$, the function $T\mapsto J_{\fo}^{T,\tG}(f)$ extends to a polynomial of degree at most $\dim\mathfrak a_{M_0}^{G}$. 
Its coefficients are continuous linear forms on $\cC(\tG(\A)^1;K_0)$. 
For each monomial in a fixed basis, the sum over $\fo$ of the absolute values of the corresponding coefficients satisfies a uniform seminorm bound with polynomial dependence on the level of $K_0$.

The twist also affects the choice of main term in the approximation of the truncated kernel. 
In the untwisted case, the only nonzero contribution for which the alternating sum over parabolic subgroups consists of a single term is the one indexed by $(G,G)$ (see Remark \ref{rem:untwistedmain}). 
In the twisted case, additional one-term contributions can occur. 
Example~\ref{ex:cga-gl3-nondecay} shows for $\GL_3$ that such a contribution need not decay as $d(T)\to\infty$, even for a fixed compactly supported test function. We therefore include all one-term contributions in the main kernel $m_{\fo}^T$ (see \S~\ref{sec:maincoarse} and \S~\ref{sec:mainrefinedclass} for the definition) and obtain the approximation estimate
\[
 \sum_{\fo\in\cO}
 \int_{G(\Q)\backslash G(\A)^1}
   |k_{\fo}^T(x,f)-m_{\fo}^T(x,f)|\,dx \le
 C\operatorname{level}(K_0)^s
 (1+\|T\|)^r e^{-d(T)}
 \|f\|_{\tG(\A)^1,m}
\]
whenever $d(T)>d_*$.
Here $k_{\fo}^T$ is the truncated kernel, $\|f\|_{\tG(\A)^1,m}$ is a seminorm defined using the $L^1$ norms of finitely many right derivatives, and $C,m,s,r$ are independent of $K_0$, $f$, $T$, and the support of $f$.
This revises the approximation estimate of \cite{Par19} and gives a version that holds uniformly for test functions
without compact support.

The proof builds on the methods of Finis--Lapid \cites{FL11semi,FL16}, with particular attention to the combinatorics of parabolic subgroups and the integration regions arising from the twist.
Estimates for weighted parabolic integrals and alternating sums over unipotent radicals reduce the main analytic problem to bounding truncated integrals.
To bound these integrals, we use the Bruhat decomposition and the Abel transform and choose an auxiliary parameter $\sigma_T$ in the Weyl polytope. 
This choice gives the required positivity in every unbounded direction of integration while preserving the exponential estimate on the bounded part of the domain, at the cost of only a polynomial factor in $1+\|T\|$. 
This construction allows us to dispense with the root cone lemma (see Remark \ref{rem:compareFL}).

On the spectral side, the author’s joint work with M\"uller \cite{MW26} gives a detailed proof of the twisted analogue of the absolute convergence theorem of Finis–Lapid–Müller. 
Together with the continuity and uniform estimates established here on the geometric side, this provides a basis for applying the twisted trace formula to test functions of noncompact support. 
The revised approximation estimate enters directly into the proof of the limit multiplicity formula in \cite{TW26}. 
Remark~\ref{rem:TW} explains how to adapt the proof of the limit multiplicity formula in \cite{TW26} by replacing Parab's estimate \cite{Par19} with the revised approximation estimate proved here.
The author and M\"uller also plan to use the revised estimate in their ongoing work on a Weyl law for self-dual cuspidal automorphic representations, following the strategy of Finis and Matz \cite{FM21}. 
For a different application of noncompactly supported test functions, Parab \cite{Par19}*{\S~8} discusses sums of residues of certain Rankin–Selberg $L$-functions in connection with Langlands’ Beyond Endoscopy proposal.

The paper is organized as follows.
In \S~\ref{sec:preliminaries}, we collect notation and analytic preliminaries.
In \S~\ref{sec:cga}, we define the main term and prove the decay estimate for test functions supported in a fixed compact set. 
In \S~\ref{gr:section}, we prove the continuity of the coarse geometric expansion, and in \S~\ref{fg:section}, we establish analogous results for its refinement using induced conjugacy classes.

\medskip
\noindent
{\bf Acknowledgment}: The author is partially supported by JSPS Grant-in-Aid for Challenging Research (Exploratory) No.26K22252.

\medskip
\noindent
{\bf AI disclosure}: During the preparation of this paper, the author used ChatGPT (OpenAI) to assist with drafting the English text, improving the exposition, and developing mathematical arguments. In particular, Example~\ref{ex:cga-gl3-nondecay}, which concerns decay estimates, was found by ChatGPT and motivated the writing of this paper. The author has checked all content generated by ChatGPT and takes full responsibility for the content of the paper.

\section{Preliminaries}
\label{sec:preliminaries}

In this section, we introduce the notation and the coarse geometric expansion of the twisted trace formula.

\subsection{Algebraic groups}

Let $\A$ be the adele ring of $\Q$, and let $\A_f$ be its finite
adele ring. For each place $v$ of $\Q$, let $|\cdot|_v$ be the
normalized absolute value on $\Q_v$.
Write $|\cdot|_{\A}=\prod_v|\cdot|_v$ for the idele norm.
We fix the following data.
\begin{itemize}
\item $G'$ is a reductive algebraic group over $\Q$.
\item $G$ is the connected component of $1$ in $G'$. It is known that $G$ is a normal closed subgroup of $G'$, of finite index, and defined over $\Q$.
\item $\tG$ is a connected component of $G'$ such that
$\tG(\Q)\neq\emptyset$.
\item Fix a minimal parabolic subgroup $P_0$ of $G$ over $\Q$.
\item Fix a Levi subgroup $M_0$ of $P_0$ over $\Q$. 
\item Define a homomorphism $\Ad$ from $G'$ to the automorphism group of $G$ by $\Ad(g)h \coloneqq ghg^{-1}$ for $g \in G'$ and $h \in G$.
\item  There exists an element $\delta\in \tG(\Q)$ such that $\Ad(\delta)(P_0)=P_0$ and $\Ad(\delta)(M_0)=M_0$.
\item In the following, we fix one such $\delta \in \tG(\Q)$ and set $\theta \coloneqq \Ad(\delta)$.
      Then $\theta$ is an automorphism of $G$ defined over $\Q$, and $P_0$ and $M_0$ are $\theta$-stable, 
      that is, $\theta(P_0)=P_0$ and $\theta(M_0)=M_0$.
\end{itemize}

Let $T_0$ be the maximal $\Q$-split torus in the center of
$M_0$. 
Thus $M_0$ is the centralizer of $T_0$ in $G$.
Let $\cL^G(M_0)$ denote the set of Levi $\Q$-subgroups containing $M_0$, that is, the finite set of centralizers of $\Q$-subtori of $T_0$. 
For $M\in\cL^G(M_0)$, we use the following notation.
\begin{itemize}        
    \item $T_M$ is the maximal $\Q$-split subtorus of the identity component of the center of $M$.
    \item $\cL^G(M)$ is the set of Levi $\Q$-subgroups containing $M$. 
    \item $\cP^G(M)$ is the set of parabolic $\Q$-subgroups with Levi part $M$. 
    \item $\cF^G(M)\coloneqq\coprod_{L\in\cL^G(M)} \cP^G(L)$ is the finite set of parabolic $\Q$-subgroups of $G$ containing $M$. 
    \item For each $P\in\cP^G(M)$, let $N_P$ be the unipotent radical of $P$, and let $M_P$ be the Levi subgroup of $P$ such that $M_P\in \cL^G(M_0)$. 
    \item A parabolic subgroup $P\in\cF^G(M_0)$ is called standard if it contains $P_0$.
\end{itemize}

We explain the parabolic subsets of $\tG$ and related notions.
\begin{itemize}
    \item For $P\in\cF^G(M_0)$, we write $\tP$ for the normalizer of $P$ in $\tG$. If $\tP$ is nonempty, we call $\tP$ a parabolic subset of $\tG$. In this case, $\tP(\Q)$ is nonempty.
    \item Set $\tM_0\coloneqq M_0\delta$ and $\cF^{\tG}(\tM_0)\coloneqq\{ \tP \mid P\in\cF^G(M_0)$ and $\tP\supset\tM_0 \}$. 
    All elements in $\cF^{\tG}(\tM_0)$ are parabolic subsets of $\tG$. 
    In addition, any $\tP\in\cF^{\tG}(\tM_0)$ satisfies $\tP=P\delta$ and $\theta(P)=P$. 
    \item For $\tP\in \cF^{\tG}(\tM_0)$, we have $\theta(M_P)=M_P$ and $\theta(N_P)=N_P$. Hence, we have a Levi decomposition $\tP = \tM_P N_P$, where $\tM_P\coloneqq M_P\, \delta$. 
    We call this $\tM_P$ a Levi subset of $\tG$.
\end{itemize}
For a Levi subset $\tM$ of $\tG$, we use the following notation. 
\begin{itemize}
    \item $T_{\tM}$ is the $\Q$-split part of the identity component of the centralizer of $\tM$ in $M$. 
    Then, $T_{\tM}$ is the identity component of $\{ t\in T_M \mid \theta(t)=t \}$.
    \item $\cF^{\tG}(\tM)\coloneqq \{ \tP\in\cF^{\tG}(\tM_0) \mid \tM_P\supset \tM\}$. 
    \item $\cL^{\tG}(\tM)\coloneqq \{ \tM_P \mid \tP\in \cF^{\tG}(\tM)\}$. 
    \item $\cP^{\tG}(\tM) \coloneqq \{ \tP\in \cF^{\tG}(\tM) \mid \tM_P=\tM \}$.
    \item A parabolic subset $\tP\in\cF^{\tG}(\tM_0)$ is standard when $P\supset P_0$.
\end{itemize}

\subsection{Real vector spaces}

For $M\in\cL^G(M_0)$, let $X^*(M)$ denote the lattice of
$\Q$-rational characters of $M$, and put
\[
  \fa_M^*\coloneqq X^*(M)\otimes_{\Z}\R,
  \qquad \fa_M\coloneqq \Hom_{\R}(\fa_M^*,\R).
\]
For $P\in\cP^G(M)$, set $\fa_P\coloneqq \fa_M$.
If $M\subset L$, use the usual decompositions
\[
  \fa_M=\fa_L\oplus\fa_M^L,
  \qquad \fa_M^*=\fa_L^*\oplus(\fa_M^L)^*,
\]
see \cite{Art05}*{\S5}. 
Here, $\fa_M^L$ (resp. $(\fa_M^L)^*$) is the annihilator of $\fa_L^*$ in $\fa_M$ (resp. $\fa_L$ in $\fa_M^*$).
For $P\subset R$, write $\fa_P^R=\fa_{M_P}^{M_R}$.

Fix a maximal compact subgroup $K=K_\infty K_f$ of $G(\A)$
that is admissible relative to $M_0$, see \cite{Art81}*{\S1}. 
Note that $\theta(K)$ is also admissible relative to $M_0$. 
Then, the Iwasawa decomposition $G(\A)=P(\A)\, K$ holds for any $P\in\cF^G(M_0)$. 
Let $A_0$ denote the connected component of $1$ in $T_0(\R)$. 
For $M\in \cL^G(M_0)$ and $P\in\cP^G(M)$, put
\[
   A_M=A_P\coloneqq A_0\cap T_M(\R).
\]
For $m\in M(\A)$, define $H_M(m)\in\fa_M$ by
\[
  \langle\chi,H_M(m)\rangle=\log|\chi(m)|_{\A},
  \qquad \chi\in X^*(M).
\]
For $P\in\cP^G(M)$, extend $H_M$ to $H_P\colon G(\A)\to\fa_P$
by the Iwasawa decomposition:
\[
  H_P(nmk)=H_M(m),
  \qquad n\in N_P(\A),\quad m\in M(\A),\quad k\in K.
\]
The map $H_M\colon A_M\to\fa_M$ is an isomorphism.

Following the notation for the spectral side, put
\begin{equation*}
  G(\A)^1=\{x\in G(\A)\mid H_G(x)=0\},
  \qquad \tG(\A)^1=G(\A)^1\delta.
  \label{eq:lw911-norm-one}
\end{equation*}
Then
\[
  G(\A)=A_G\times G(\A)^1,
  \qquad A_G\backslash G(\A)\simeq G(\A)^1.
\]
The product formula implies $G(\Q)\subset G(\A)^1$ and
$\tG(\Q)\subset\tG(\A)^1$.

Let $\tM\in\cL^{\tG}(\tM_0)$ and $\tP\in \cP^{\tG}(\tM)$.
For $\chi\in X^*(T_M)$, we set $\theta(\chi)\coloneqq\chi\circ\theta^{-1}$. 
This defines an action of $\theta$ on $\fa_M^*$ and, by duality, also an action on $\fa_M$.
Its action on $\fa_M$ has finite order.
Indeed, some positive power $\delta^\ell$ belongs to $G$;
since it normalizes both $P$ and $M$, it belongs to $M$.
This does not require $\theta$ itself to have finite order on $G$.
We set
\[
  \fa_{\tM}\coloneqq \fa_M^\theta=\{H\in \fa_M \mid \theta(H)=H \},\quad 
  \fa_{\tM}^*\coloneqq \fa_M^{*,\theta}=\{\lambda\in \fa_M^* \mid \theta(\lambda)=\lambda \}.
\]
For $\tL\in\cL^{\tG}(\tM)$, set
\[
  \fa_{\tM}^{\tL}\coloneqq (\fa_M^L)^\theta=\{H\in \fa_M^L \mid \theta(H)=H \},\quad 
  (\fa_{\tM}^{\tL})^*\coloneqq (\fa_M^L)^{*,\theta}=\{\lambda\in (\fa_M^L)^* \mid \theta(\lambda)=\lambda \}.
\]
For $\tP\subset\tQ$, put $\fa_{\tP}^{\tQ}\coloneqq \fa_{\tM_P}^{\tM_Q}$ and $(\fa_{\tP}^{\tQ})^*\coloneqq (\fa_{\tM_P}^{\tM_Q})^*$.
Then, we have
\[
\fa_{\tP}=\fa_{\tQ}\oplus \fa_{\tP}^{\tQ}, \quad \fa_{\tP}^*=\fa_{\tQ}^*\oplus (\fa_{\tP}^{\tQ})^*.
\]
Set $\fa_M^{1-\theta}\coloneqq \langle H-\theta(H)\mid H\in\fa_M\rangle$, and then 
\begin{equation*}
  \fa_M=\fa_{\tM}\oplus \fa_M^{1-\theta}.
  \label{eq:lw911-central-assumption}
\end{equation*}
Note that $\fa_{\tM}^*$ is the annihilator of $\fa_M^{1-\theta}$ in $\fa_M^*$, and $(\fa_M^{1-\theta})^*$ is identified with the annihilator of $\fa_{\tM}$ in $\fa_M^*$. 
Hence, $\fa_M^*=\fa_{\tM}^*\oplus (\fa_M^{1-\theta})^*$.

\subsection{On the normalization of measures}\label{sec:measure}


We retain the groups $K$, $A_0$, and $A_M$ fixed above.
For $P\in\cP^G(M)$, let $\delta_P$ be the modular character
of $P(\A)$, and let $W_0^G$ be the Weyl group of $(G,T_0)$.

Let us fix an inner product on $\fa_{M_0}$, which is $W_0^G$-invariant and $\theta$-invariant.
Then, for the orthogonal complement $(\fa_{\tM})^\perp$ of $\fa_{\tM}$ in $\fa_M$, we have
\[
\fa_M^{1-\theta}=(\fa_{\tM})^\perp.
\]
This inner product endows $\fa_{M_0}$ with the structure of a Euclidean space. 
In addition, it determines Haar measures on $\fa_{M}^L$ and $\fa_{\tM}^{\tL}$.
The dual measures on $(\fa_{M}^{L})^*$ and $(\fa_{\tM}^{\tL})^*$ are then fixed accordingly.
Moreover, via the isomorphism $H_M\colon A_M\to \fa_M$, a measure on $A_M$ is also fixed.

Let us fix Haar measures on groups following \cite{Art78}*{\S1}. 
First, fix a Haar measure $dx$ on $G(\A)$. 
Normalize the Haar measure $dk$ on $K$ so that the volume of $K$ is $1$. 
For each $P\in\cF^G(M_0)$, normalize the Haar measure $dn$ on $N_P(\A)$ so that the volume of $N_P(\Q)\bs N_P(\A)$ is $1$. 
Note that the Iwasawa decomposition $G(\A)=N_P(\A) M_P(\A)K$ holds. 
Then normalize the Haar measure $dm$ on $M_P(\A)$ so that
\[
\int_{G(\A)} f(x)\, dx =\int_{N_P(\A)} \int_{M_P(\A)} \int_K f(nmk) \, \delta_P(m)^{-1}  dk\, dm \, dn, \qquad  f\in C_c(G(\A)) 
\]
holds.
Induce a measure $d\widetilde{x}$ on $\tG(\A)$ from the measure $dx$ on $G(\A)$ by
\[
\int_{\tG(\A)}f(\widetilde{x})\, d\widetilde{x}=\int_{G(\A)} f(x\delta) \, dx,\qquad f\in C_c(\tG(\A)).
\]

Let $M\in\cL^G(M_0)$ and $P\in \cP^G(M)$. 
Set
\[
A_0^G\coloneqq\{a\in A_0 \mid H_G(a)=0\}, \quad A_M^G\coloneqq A_0^G\cap A_M. 
\]
Let $\exp\colon\fa_{M_0}\to A_0$ be the inverse of $H_{M_0}$. In particular, $\exp(\fa_{M}^G)=A_M^G$. 
For $P\in\cF^G(M_0)$, $\rho_P$ denotes half the sum, with multiplicities, of the positive
roots corresponding to $(T_M,N_P)$.
Then, 
\[
\delta_P(a)=e^{2\rho_P(H_M(a))} ,\quad a\in A_M.
\]
We will often use
\[
\int_{G(\A)^1} f(x)\, dx =\int_{N_P(\A)} \int_{M(\A)^1}\int_{A_M^G} \int_K f(nmak) \, e^{-2\rho_P(H_P(a))}  dk\, da\, dm \, dn, \quad  f\in C_c(G(\A)^1) .
\]

\subsection{Roots and weights}

Take a standard parabolic subgroup $P\in\cF^G(M_0)$ with $P\supset P_0$, and let $M=M_P$. 
For the definitions of the following notation, see \cite{Art05}*{\S5}. 
\begin{itemize}
    \item $\Delta_P$ is the set of simple roots of $T_M$ in $N_P$. 
    \item For $\alpha\in\Delta_P$, let $\alpha^\vee$ denote the coroot corresponding to $\alpha$ in the sense of \cite{Art05}*{\S5}. Set $\Delta_P^\vee\coloneqq \{ \alpha^\vee \mid \alpha\in\Delta_P\}$.     
    \item When $R\in\cF^G(M)$ and $R\supset P$, we set
\[
  \Delta_P^R\coloneqq \{\alpha\in\Delta_P\mid \alpha|_{\fa_R}=0\},\quad 
  (\Delta_P^R)^\vee\coloneqq \{\alpha^\vee \mid \alpha\in \Delta_P^R\}.
\]
    \item $\widehat\Delta_P^R$ is the basis of $(\fa_P^R)^*$ dual to $(\Delta_P^R)^\vee$.
    \item $(\widehat{\Delta}_P^R)^\vee$ is the basis of $\fa_{P}^{R}$ dual to $\Delta_{P}^{R}$.
\end{itemize}
We regard these roots and weights as forms on $\fa_{M_0}$ by
composing with the projection onto $\fa_P^R$.
For simplicity, we set $\widehat\Delta_P\coloneqq \widehat\Delta_P^G$. 
The map $P\mapsto \Delta_{P_0}^P$ gives a bijection between $\{P\in\cF^G(M_0)\mid P\supset P_0\}$ and the subsets of $\Delta_{P_0}$.

Suppose that $P$ is $\theta$-stable.
$\theta$ acts bijectively on $\Delta_P$,
since $N_P$ is $\theta$-stable.
For each $\alpha\in\Delta_P$, let $\ell_\alpha$ denote
the cardinality of its $\theta$-orbit. We set
\[
  \talpha \coloneqq
    \frac{1}{\ell_\alpha}
    \sum_{j=0}^{\ell_\alpha-1}\theta^j(\alpha),
  \qquad
  \talpha^\vee \coloneqq
    \sum_{j=0}^{\ell_\alpha-1}\theta^j(\alpha^\vee).
\]
Put
\[
  \Delta_{\tP}
     \coloneqq \{\widetilde\alpha \mid \alpha\in\Delta_P\},
  \quad
  \Delta_{\tP}^\vee
     \coloneqq \{\widetilde\alpha^\vee\mid \alpha\in\Delta_P\}.
\]

Let $\widehat\Delta_{\tP}$ (resp. $\widehat\Delta_{\tP}^\vee$) be the basis of $(\fa_{\tP}^{\tG})^*$ (resp. $\fa_{\tP}^{\tG}$) dual to $\Delta_{\tP}^\vee$ (resp. $\Delta_{\tP}$). 
Then, similarly to $\Delta_{\tP}$, the set $\widehat\Delta_{\tP}$ consists of the $\theta$-orbit averages of the elements of $\widehat{\Delta}_P$.
\begin{remark}
In \cite{MW26}, $\widetilde\alpha^\vee$ was defined as the $\theta$-orbit average of $\alpha^\vee$.
Here we instead use the $\theta$-orbit sum. 
Consequently, the basis dual to $\Delta_{\tP}^\vee$ consists of the $\theta$-orbit averages of the elements
of $\widehat\Delta_P$. 
The elements of $\Delta_{\tP}$ and $\widehat\Delta_{\tP}$
are therefore the images of the elements of
$\Delta_P$ and $\widehat\Delta_P$, respectively,
under the natural projection $\fa_P^* \to\fa_{\tP}^*$.
This agrees with the normalization of roots and weights in \cite{LW13}*{\S2.7}.
\end{remark}

For $\tQ\in\cP^{\tG}(\tL)$, $Q\supset P$, we set
    \[
    \Delta_{\tP}^{\tQ}\coloneqq \{\talpha\in\Delta_{\tP} \mid \talpha|_{\fa_{\tQ}}=0 \}, \quad (\Delta_{\tP}^{\tQ})^\vee\coloneqq \{\talpha^\vee  \mid \talpha\in \Delta_{\tP}^{\tQ} \}.
    \]
Let $\widehat\Delta_{\tP}^{\tQ}$ (resp. $(\widehat\Delta_{\tP}^{\tQ})^\vee$) be the basis of
$(\fa_{\tP}^{\tQ})^*$ (resp. $\fa_{\tP}^{\tQ}$) dual to $(\Delta_{\tP}^{\tQ})^\vee$ (resp. $\Delta_{\tP}^{\tQ}$).

We record the following two basic facts, cf. \cite{LW13}*{Lemmes 2.7.1, 2.7.2}.
\begin{itemize}
    \item Fix $\tP \in \cF^{\tG}(\tM_0)$.
    The map $\tQ \longmapsto \Delta_{\tP}^{\tQ}$ gives a bijection from the set $\{\, \tQ \in \cF^{\tG}(\tM_0) \mid \tQ \supset \tP \,\}$ onto the set of subsets of $\Delta_{\tP}$.
    \item Let $\tP, \tQ \in \cF^{\tG}(\tM_0)$ and assume $\tP \subset \tQ$. Then $\Delta_{\tP}^{\tQ}$ (resp. $(\Delta_{\tP}^{\tQ})^\vee$) forms a basis of $(\fa_{\tP}^{\tQ})^*$ (resp. $\fa_{\tP}^{\tQ}$).
\end{itemize}

We introduce two types of functions on $\fa_{M_0}$.
\begin{itemize}
    \item $\tau_P^R$ is the characteristic function of $\{H\in\fa_{M_0}\mid  \alpha(H)>0, \;\; \forall \alpha\in\Delta_P^R\}$ on $\fa_{M_0}$.  
    \item $\widehat\tau_P^R$ is the characteristic function of $\{H\in\fa_{M_0}\mid  \varpi(H)>0, \;\; \forall \varpi\in\widehat\Delta_P^R\}$ on $\fa_{M_0}$.  
    \item $\tau_{\tP}^{\tR}$ is the characteristic function of $\{H\in\fa_{M_0}\mid  \alpha(H)>0, \;\; \forall \alpha\in\Delta_{\tP}^{\tR}\}$ on $\fa_{M_0}$.  
    \item $\widehat\tau_{\tP}^{\tR}$ is the characteristic function of $\{H\in\fa_{M_0}\mid  \varpi(H)>0, \;\; \forall \varpi\in\widehat\Delta_{\tP}^{\tR}\}$ on $\fa_{M_0}$.  
\end{itemize}

\subsection{The function
             \texorpdfstring{$\widetilde\sigma_Q^R$}{sigma}}

Let $\ell$ be the order of the action of $\theta$ on $\fa_{M_0}^*$, hence we assume that $\theta^\ell$ acts trivially on $\fa_{M_0}^*$.
Let $Q$, $R\in\cF^G(M_0)$ with $Q\subset R$.
Suppose that $Q$ and $R$ are standard.
Define standard parabolic subgroups $Q^+$ and $R^-$ by
\[
  \Delta_{P_0}^{Q^+}=\bigcup_{j=0}^{\ell-1}\theta^j(\Delta_{P_0}^{Q}),
    \quad 
  \Delta_{P_0}^{R^-}=\bigcap_{j=0}^{\ell-1}\theta^j(\Delta_{P_0}^{R}).
\]
Thus $Q^+$ is the smallest $\theta$-stable standard parabolic
subgroup containing $Q$, and $R^-$ is the largest such subgroup
contained in $R$. 
There exists a standard parabolic subset $\tP\in\cF^{\tG}(\tM_0)$ with
$Q\subset P\subset R$ if and only if $Q^+\subset R^-$.

Suppose that $Q^+\subset R^-$, and choose any
$\tP$ in this interval.
Define $\widetilde\sigma_Q^R$ to be the characteristic
function of the set of $H\in\fa_{M_0}$ satisfying
\[
\begin{cases}
  \alpha(H)>0
    &(\forall\alpha\in\Delta_Q^R),\\
  \alpha(H)\leq0
    &(\forall\alpha\in\Delta_Q\setminus\Delta_Q^R),\\
  \tvarpi(H)>0
    &(\forall\tvarpi\in\widehat\Delta_{\tP}).
\end{cases}
\]
By \cite{LW13}*{Lemma~2.11.3}, this definition is independent of the choice of $\tP$.
If $Q^+\not\subset R^-$, set $\widetilde\sigma_Q^R\coloneqq 0$.

The following identities will be useful:
\[
  \widetilde\sigma_Q^G=\tau_Q^G,
  \qquad \widetilde\sigma_Q^Q=0\quad(Q\neq G),
  \qquad \widetilde\sigma_G^G=1.
\]
Moreover, if $Q\subset P$ and $\tP\in\cF^{\tG}(\tM_0)$, then
\begin{equation}
  \sum_{\substack{R\in\cF^G(M_0) \\ R\supset P}}\widetilde\sigma_Q^R
      =\tau_Q^P\widehat\tau_{\tP}.
  \label{eq:lw911-sigma-identity}
\end{equation}
See \cite{LW13}*{Lemma~2.11.5}.

\subsection{Siegel sets and partition formula}

Suppose $\Delta_{P_0}\neq \emptyset$. 
For $T\in \fa_{M_0}^G$, we set
\[
d(T)\coloneqq \min_{\alpha\in\Delta_{P_0}} \alpha(T).
\]
We fix a parameter $T_1\in\fa_{M_0}^G$ such that $d(-T_1)$ is sufficiently large. 
When $\Delta_{P_0}=\emptyset$, we put $d(T)\coloneqq+\inf$.
Following \cite{FL16}*{\S2.1}, we introduce some formulas.

Write $H_0\coloneqq H_{P_0}$.
Take an arbitrary $P\in\cF^G(M_0)$ containing $P_0$.
Set 
\[
\fS_{T_1}^P\coloneqq\{x\in G(\A)\mid\tau_{P_0}^P(H_0(x)-T_1)=1\}. 
\]
By reduction theory, $P(\Q)\fS_{T_1}^P=G(\A)$.

We define $F^P(x,T)$ to be the characteristic function in $x$ of
    \[
  \left\{x\in G(\A)\ \middle|\
      \varpi\bigl(H_0(\delta x)-T\bigr)\leq 0
      \text{ for every }\delta\in P(\Q) \;\; 
      \text{and every }\varpi\in\widehat\Delta_{P_0}^P
  \right\}.  
    \]
Note that $F^P(\cdot,T)$ is left $P(\Q)$-invariant and right
$K$-invariant, and $F^{P_0}(x,T)=1$.  
\begin{remark}
\label{rem:truncation-Arthur-compatibility}
Fix a compact subset $\omega\subset N_{P_0}(\A)M_0(\A)^1$ such that
the Siegel sets
\[
  \cS^P(T_1,\omega)
  \coloneqq
  \left\{pak\ \middle|\
      p\in\omega,\ a\in A_0,\ k\in K,\;\; 
      \alpha(H_0(a)-T_1)>0
      \text{ for all }\alpha\in\Delta_{P_0}^P
  \right\}
\]
($P\in\cF^G(M_0)$ with $P\supset P_0$) satisfy $P(\Q)\cS^P(T_1,\omega)=G(\A)$.
Arthur defines $F^P(\cdot,T)$ using the image in
$P(\Q)\backslash G(\A)$ of
\begin{equation*}
  \cS^P(T_1,T,\omega)
  \coloneqq
  \left\{y\in\cS^P(T_1,\omega)\ \middle|\
    \varpi(H_0(y)-T)\leq0
    \text{ for all }\varpi\in\widehat\Delta_{P_0}^P
  \right\}.
\end{equation*}
See \cite{Art78}*{\S6, p.~941} and \cite{Art05}*{\S8}.
There exists $d_0>0$ such that, whenever
$d(T)>d_0$, for any $P\in\cF^G(M_0)$ containing $P_0$ and any $y\in \fS_{T_1}^P$, we have
\[
\forall \varpi\in\widehat\Delta_{P_0}^P, \ \  \varpi(H_0(y)-T)\leq0  \quad\Longleftrightarrow
  \quad \forall \xi\in P(\Q),\ \  \forall\varpi\in\widehat\Delta_{P_0}^P, \ \ \varpi(H_0(\xi y)-T)\leq 0.
\]
This equivalence allows us to show that Arthur's definition of $F^P(\cdot,T)$ agrees with the above definition of Finis--Lapid for $d(T)>d_0$.
\end{remark}

Increasing $d_0$ if necessary, we obtain the following partition formula. 
\begin{lemma}
\label{lem:partitionformula}
Take a standard parabolic subgroup $P\in\cF^G(M_0)$. 
For every $T\in\fa_{M_0}^G$ satisfying $d(T)>d_0$ and every
$x\in G(\A)$, we have
\begin{equation*}
  \sum_{\substack{Q\in\cF^G(M_0)\\ P_0 \subset Q\subset P}}
  \ \sum_{\xi\in Q(\Q)\backslash P(\Q)}
  F^Q(\xi x,T)\,
  \tau_Q^P\bigl(H_Q(\xi x)-T\bigr)
  =1.
\end{equation*}
\end{lemma}
\begin{proof}
See \cite{Art78}*{Lemma~6.4}.
\end{proof}


\subsection{Kernel functions}

Let $f\in C_c^\infty(\tG(\A)^1)$.
For $\tP\in\cF^{\tG}(\tM_0)$ containing $\tP_0$, we set
\begin{equation*}
  K_{\tP}(x,y)
    \coloneqq \sum_{\gamma\in\tM_P(\Q)}
       \int_{N_P(\A)} f(x^{-1}\gamma n y)\,dn, \quad x, y\in G(\A)^1,
\end{equation*}
\[
  k^{T}(x,f)
   \coloneqq\sum_{\substack{\tP\in\cF^{\tG}(\tM_0)\\P\supset P_0}}
      (-1)^{\dim\fa_{\tP}^{\tG}}
      \sum_{\xi\in P(\Q)\backslash G(\Q)}
       K_{\tP}(\xi x,\xi x)\, \widehat{\tau}_{\tP}(H_P(\xi x)-T).
\]

Every $\gamma\in\tG(\Q)$ has its usual Jordan decomposition
in the algebraic group $G'$:
\[
  \gamma=\gamma_s\gamma_u=\gamma_u\gamma_s,
  \qquad \gamma_s\in\tG(\Q),\quad \gamma_u\in G(\Q),
\]
where $\gamma_s$ is semisimple and $\gamma_u$ is unipotent.
Two elements of $\tG(\Q)$ belong to the same coarse class
if their semisimple parts are conjugate under $G(\Q)$.
Let $\tcO$ denote the set of coarse classes.
For $\tfo\in\tcO$ and $\tP\in\cF^{\tG}(\tM_0)$ containing $\tP_0$, we set
\begin{equation*}
  K_{\tP,\tfo}(x,y)
    \coloneqq \sum_{\gamma\in\tM_P(\Q)\cap\tfo}
       \int_{N_P(\A)} f(x^{-1}\gamma n y)\,dn, \quad x, y\in G(\A)^1,
\end{equation*}
\begin{equation*}
  k_{\tfo}^{T}(x,f)
   \coloneqq\sum_{\substack{\tP\in\cF^{\tG}(\tM_0)\\P\supset P_0}}
      (-1)^{\dim\fa_{\tP}^{\tG}}
      \sum_{\xi\in P(\Q)\backslash G(\Q)}
       K_{\tP,\tfo}(\xi x,\xi x)\, \widehat{\tau}_{\tP}(H_P(\xi x)-T).
\end{equation*}
Then, $K_{\tP}(x,y)=\sum_{\tfo\in\tcO}K_{\tP,\tfo}(x,y)$ and $k^{T}(x,f)=\sum_{\tfo\in\tcO}k^{T}_{\tfo}(x,f)$.

\begin{proposition}
\label{prop:lw911}
Fix a compact subset $\Omega\subset\tG(\A)^1$.
There is a constant $c_\Omega$, depending on $\Omega$ and the
fixed reduction data, with the following property.
Let $Q$, $R\in\cF^G(M_0)$ and suppose $P_0\subset Q\subset R$. 
Let $T\in\fa_{M_0}^G$, $\tfo\in\tcO$ and $f\in C_c^\infty(\tG(\A)^1)$ with $\supp(f)\subset\Omega$.
If $d(T)\geq c_\Omega$, then the integral
\[
\sum_{\tfo\in\tcO}\int_{Q(\Q)\backslash G(\A)^1}
       F^Q(x,T)\,
        \widetilde\sigma_Q^R(H_0(x)-T)\, 
       \Bigl|
         \sum_{\substack{\tP\in\cF^{\tG}(\tM_0) \\ \widetilde Q^+\subset\tP
                                 \subset\widetilde R^-}}
           (-1)^{\dim\fa_{\tP}^{\tG}}
           K_{\tP,\tfo}(x,x)
       \Bigr|\,dx
\]
is finite. 
Note that, in the above sum, the integral vanishes for all but finitely many $\tfo\in\tcO$.
\end{proposition}
\begin{proof}
See \cite{LW13}*{Proposition~9.1.1}.
\end{proof}

Let $Q$, $R\in\cF^G(M_0)$ such that $P_0\subset Q\subset R$. 
When $Q^+\subset R^-$, we set
\begin{equation*}
 B_{\tfo,Q,R}^T(x,f)
 \coloneqq \sum_{\xi\in Q(\Q)\backslash G(\Q)}
    F^Q(\xi x,T)\tsigma_Q^R(H_0(\xi x)-T)
    \sum_{\substack{\tP\in\cF^{\tG}(\tM_0) \\ \tQ^+\subset \tP\subset \tR^-}}
       (-1)^{\dim\fa_{\tP}^{\tG}}\, K_{\tP,\tfo}(\xi x,\xi x).
\end{equation*}
When $Q^+\not\subset R^-$, we set $B_{\tfo,Q,R}^T(x,f)\coloneqq 0$. 
Applying Lemma~\ref{lem:partitionformula} to $P$, $Q$ and using \eqref{eq:lw911-sigma-identity} gives
\begin{equation*}
 k_{\tfo}^T(x,f)=\sum_{Q\subset R}B_{\tfo,Q,R}^T(x,f).
\end{equation*}
Hence, by Proposition \ref{prop:lw911}, we obtain
\[
  \sum_{\tfo\in\tcO}\int_{G(\Q)\backslash G(\A)^1}|k_{\tfo}^T(x,f)|\,dx    <\infty.
\]
Thus the integrals
\[
  J_{\tfo}^{T,\tG}(f)
     \coloneqq \int_{G(\Q)\backslash G(\A)^1}
          k_{\tfo}^T(x,f)\,dx, \qquad   J^{T,\tG}(f)
     \coloneqq \int_{G(\Q)\backslash G(\A)^1}
          k^T(x,f)\,dx
\]
are well defined for sufficiently regular $T$, and then $J^{T,\tG}(f)=\sum_{\tfo\in\tcO}J_{\tfo}^{T,\tG}(f)$.

\subsection{Fr\'echet spaces}\label{sec:Frechet}

We recall the spaces of test functions used in
\cite{FL16}*{\S2.2}.
Put
\[
 G(\R)^1=G(\R)\cap G(\A)^1,\qquad
 \fg^1=\Lie(G(\R)^1),\qquad
 \fg_\C^1=\fg^1\otimes_\R\C.
\]
For $X\in\fg^1$, let
\[
 (\phi*X)(x)
   =\left.\frac{d}{dt}\phi(x\exp(tX))\right|_{t=0},
 \qquad x\in \tG(\A)^1.
\]
This action extends to the universal enveloping algebra
$\cU(\fg_\C^1)$.

Fix an ordered basis of $\fg^1$.
Let $X_1,\ldots,X_d$ be the fixed ordered $\C$-basis of $\fg_\C^1$.
The standard filtration on $\cU(\fg_\C^1)$ is defined by
\[
 \cU(\fg_\C^1)_{\leq m}
 \coloneqq
 \operatorname{span}_\C
 \left\{
 Y_1\cdots Y_r\mid
 0\leq r\leq m,\quad Y_1,\ldots,Y_r\in\fg_\C^1
 \right\},
\]
where the empty product is the identity element $1$.
For a multi-index $\alpha=(\alpha_1,\ldots,\alpha_d)\in\mathbb Z_{\geq0}^d$, put $|\alpha|\coloneqq \alpha_1+\cdots+\alpha_d$ and $X^\alpha\coloneqq X_1^{\alpha_1}\cdots X_d^{\alpha_d}$ with the factors taken in the indicated order.
By the Poincar\'e--Birkhoff--Witt theorem, for each integer $m\geq0$,
\[
 \mathcal B_m
 \coloneqq
 \left\{
 X^\alpha\mid
 \alpha\in\mathbb Z_{\geq0}^d,\quad |\alpha|\leq m
 \right\}
\]
is a $\C$-basis of $\cU(\fg_\C^1)_{\leq m}$.
In particular, $\mathcal B_0=\{1\}$ and
$\mathcal B_m\subset\mathcal B_{m+1}$.
Set
\begin{equation}
 \|\phi\|_{\tG(\A)^1,m}
   =\sum_{X\in\mathcal B_m}
        \|\phi*X\|_{L^1(\tG(\A)^1)}.
 \label{eq:frechet-test-norms}
\end{equation}
For an open compact subgroup $K_0\subset G(\A_f)$, define
\[
 \cC(\tG(\A)^1;K_0)
   \coloneqq \left\{\phi\in C^\infty(\tG(\A)^1) \; \middle| \; 
       \begin{array}{l}
        \phi(xk)=\phi(x)\quad
            (x\in \tG(\A)^1,\ k\in K_0),\\
        \|\phi\|_{\tG(\A)^1,m}<\infty
            \quad\text{for every }m\geq0
       \end{array}
     \right\}.
\]
Here smoothness means smoothness on the real Lie group
components of $\tG(\A)^1/K_0$.
The increasing family of norms
\eqref{eq:frechet-test-norms} defines a Fr\'echet topology on
$\cC(\tG(\A)^1;K_0)$, independent of the chosen basis.
In particular, this space is complete, metrizable, and locally
convex.

A seminorm $\mu$ on $\cC(\tG(\A)^1;K_0)$ is continuous if and
only if there exist an integer $m\geq0$ and a constant $C>0$
such that
\[
 \mu(\phi)\leq C\|\phi\|_{\tG(\A)^1,m}
 \qquad\bigl(\phi\in\cC(\tG(\A)^1;K_0)\bigr).
\]
The order $m$ and the constant $C$ may depend on $\mu$ and
$K_0$; any uniformity in $K_0$ will be stated explicitly.
The subspace $C_c^\infty(\tG(\A)^1;K_0)$ is dense in $\cC(\tG(\A)^1;K_0)$.
As in \cite{FL16}, we also put
\[
 \cC(\tG(\A)^1)
   \coloneqq\bigcup_{K_0\subset G(\A_f)}
       \cC(\tG(\A)^1;K_0),
\]
where $K_0$ ranges over the open compact subgroups, and equip
this union with the locally convex inductive limit topology.

For a function $f$ on $\tG(\A)^1$, put
\[
 f_\delta(g)\coloneqq f(g\delta),\qquad g\in G(\A)^1.
\]
Then, we have the topological isomorphism $\cC(\tG(\A)^1;K_0)\ni f\mapsto f_\delta\in \cC(G(\A)^1; \theta(K_0))$. 

Finally, let $\Omega\subset\tG(\A)^1$ be compact.
Write $C^\infty_\Omega(\tG(\A)^1;K_0)$ for the subspace of $\cC(\tG(\A)^1;K_0)$ consisting of functions
supported in $\Omega$. 
It is closed and therefore Fr\'echet.
Its induced topology agrees with the usual topology defined by
\[
 q_m(f)\coloneqq \sum_{X\in\mathcal B_m}
             \|f*X\|_\infty,\qquad m\geq0.
\]
Indeed, 
\[
 \|f\|_{\tG(\A)^1,m}
       \leq\vol(\Omega)\,q_m(f).
\]
Conversely, local Sobolev inequalities imply that, for every
$m$, there exist $m'\geq m$ and $C_{\Omega,K_0,m}>0$ such that
\[
 q_m(f)\leq
 C_{\Omega,K_0,m}\|f\|_{\tG(\A)^1,m'}
 \qquad\bigl(\supp(f)\subset\Omega\bigr).
\]

\subsection{Levels of open compact subgroups}\label{sec:level}

Fix a faithful $\Q$-rational representation $\iota\colon G'\to\GL(E)$ and a full $\Z$-lattice
$\Lambda\subset E$, and put $\widehat\Lambda=\Lambda\otimes_\Z\widehat\Z$.
For $n\geq1$, define the principal congruence subgroup
\begin{equation*}
 K(n)\coloneqq \{g\in G(\A_f)\mid \iota(g)\in\GL(\widehat\Lambda),\quad
             (\iota(g)-1)\widehat\Lambda\subset n\widehat\Lambda\}.
\end{equation*}
Here only the restriction of $\iota$ to $G$ enters the definition;
using a representation of $G'$ also allows us to control
conjugation by the fixed element $\delta$ below.
The groups $K(n)$ form a neighborhood basis of the identity in
$G(\A_f)$. 
For every open compact subgroup $K_0\subset G(\A_f)$, we set
\begin{equation*}
 \level(K_0)
   \coloneqq\min\{n\geq1\mid K(n)\subset K_0\}
\end{equation*}
following \cite{FL16}*{\S2.3}.
This definition is relative to the fixed representation and
lattice; $K_0$ need not be a principal congruence subgroup or a
subgroup of $K(1)$.

We record the following lemma, which will be needed for the level estimates. Since it can be proved by an elementary argument, we omit the proof.
\begin{lemma}\label{lem:levelbound}
Let $D_f$ be a compact subset of $G(\A_f)$. 
Then, there exists a natural number $q(D_f)$ such that $x^{-1} K(r\, q(D_f)) x\subset K(r)$ holds for every $x\in D_f$ and every $r\geq1$.
\end{lemma}

\subsection{Some lemmas for semi-norms}

In this subsection, we record some lemmas for semi-norms on $\cC(\tG(\A)^1;K_0)$, which are variants of \cite{FL16}*{Lemma 2.1}.
For the proof, see \cite{FL11semi}*{\S3}. 

\begin{lemma}
    
\begin{itemize}\label{lem:1}
    \item[(1)] There exists a constant $s>0$ such that, for any $K_0$ and any $f\in\cC(\tG(\A)^1;K_0)$, we have
\begin{equation*}
 \sum_{\gamma\in \tG(\Q)}|f(\gamma)|
 \ll_{G}
 \level(K_0)^s \, \|f\|_{\tG(\A)^1, \dim G-\mathrm{rank}X^*(G)} .   
\end{equation*}

\item[(2)]
Let $D$ be a compact set of $G(\mathbb A)^1$. 
Then, we have 
\begin{equation}\label{eq:2}
 \|f(x^{-1}\cdot x)\|_{\tG(\A)^1,k} = \|f(\cdot x)\|_{\tG(\A)^1,k} \ll_{D,k} \|f\|_{\tG(\A)^1,k}    
\end{equation}
for all $f\in\cC(\tG(\A)^1)$, $k\geq0$, $x\in D$.

\item[(3)]
There exists a constant $s>0$ such that, for any $K_0$ and any $f\in\cC(\tG(\A)^1;K_0)$, there exists a right-$K$-invariant function $\tf\in\cC(\tG(\A)^1;K_f)$ satisfying $\tf(x)\geq |f(x)|$ for all $x\in\tG(\A)^1$ and
\begin{equation}\label{eq:3}
 \|\tilde f*X\|_{L^1(\tG(\mathbb A)^1)}
 \ll_X
  \level(K_0)^s\, \|f\|_{\tG(\A)^1, \dim G-\mathrm{rank}X^*(G)}     
\end{equation}
for every $X\in\mathcal U(\mathfrak g^1_\C)$.
If $f$ is compactly supported, $\tilde f$ can be chosen compactly supported.

\end{itemize}
\end{lemma}

\subsection{Alternating sums on unipotent radicals}

We record a version of \cite{FL16}*{Corollary~4.7} in which
the alternating sum is taken over $\theta$-stable parabolic
subgroups.

Let $Q\subset R$ be standard parabolic subgroups such that
$Q^+\subset R^-$.
Take $\tS\in \cF^{\tG}(\tM_0)$ with $Q^+\subset S\subset R^-$. 
The averaging map $\alpha\mapsto\talpha$ identifies the
$\theta$-orbits in
$\Delta_{P_0}^{R^-}\setminus\Delta_{P_0}^{S}$ with
$\Delta_{\tP_0}^{\tR^-}\setminus\Delta_{\tP_0}^{\tS}$.
Thus, for each
$\tbeta\in\Delta_{\tP_0}^{\tR^-}\setminus\Delta_{\tP_0}^{\tS}$,
the roots $\alpha\in\Delta_{P_0}$ satisfying
$\talpha=\tbeta$ form the corresponding $\theta$-orbit.
For $\alpha\in\Delta_{P_0}$, write
$\alpha_Q\coloneqq \alpha\circ\operatorname{pr}_{\fa_Q^G}$.
Define
\begin{equation}
 \Gamma_{Q,R;S}(T,H)
   \coloneqq\sum_{\tbeta\in
       \Delta_{\tP_0}^{\tR^-}\setminus\Delta_{\tP_0}^{\tS}}
       \ \min_{\substack{\alpha\in\Delta_{P_0}\\
                         \talpha=\tbeta}}
       \bigl\{\alpha(T)+\alpha_{Q}(H-T)\bigr\}.
 \label{eq:tfl47-gamma}
\end{equation}
The empty sum is understood to be zero.

For $S\subset P\subset R^-$, put
$N_{S}^{M_P}\coloneqq N_{S}\cap M_P$.
Take a subset $\tO$ of $\tM_{S}(\Q)$ such that $\gamma^{-1}\widetilde{O}\gamma=\widetilde{O}$ holds for any $\gamma\in M_{S}(\Q)\cap Q(\Q)$. 
For $g\in G(\A)^1$ and
$f\in\cC(\tG(\A)^1;K_0)$, consider
\begin{equation}
 \begin{split}
 \mathcal E_{Q,R;S}(g,\tO,f)
   \coloneqq {}&
   \sum_{\tS\subset \tP\subset \tR^-}
       (-1)^{\dim\fa_{\tP}^{\tG}}
       \sum_{\eta\in \tO}\sum_{\nu\in N_{S}^{M_P}(\Q)}
       \int_{N_P(\A)}
          f(g^{-1}\eta\nu ng)\,dn .
 \end{split}
 \label{eq:tfl47-alternating-sum}
\end{equation}
Here, the sum over $\tS\subset \tP\subset \tR^-$ is taken over all $\theta$-stable standard parabolic subgroups $P$ between $S$ and $R^-$.
\begin{proposition}
\label{prop:tfl47}
Let $Q$, $R \in \cF^G(M_0)$, $\tS\in\cF^{\tG}(\tM_0)$ such that $P_0\subset Q\subset R$ and $Q^+\subset S\subset R^-$.   
There exist an integer $s\geq0$ and finitely many
$X_1,\ldots,X_m\in\cU(\fg_\C^1)$ such that
\[
 |\mathcal E_{Q,R;S}(g,\tO,f)|
 \ll
 \level(K_0)^s
 e^{-\Gamma_{Q,R;S}(T,H_0(g))}
 \sum_{\eta\in\tO}\sum_{j=1}^m\int_{N_{S}(\A)}
       |(f*X_j)(g^{-1}\eta n g)|\,dn.
\]
whenever $d(T)>d_0$ and $F^Q(g,T)\, \tau_Q^{R^-}(H_Q(g)-T)=1$.
The differential operators, the exponent $s$, and the implied
constant are independent of $\tO$, $K_0$, $f$, $T$, and $g$.
Note that $F^Q(g,T)\,\widetilde\sigma_Q^R(H_0(g)-T)=1$ implies $\tau_Q^{R^-}(H_Q(g)-T)=1$.
\end{proposition}
\begin{proof}
This proposition is a twisted analogue of \cite{FL16}*{Corollary~4.7} and can be proved by the same argument as in the proof of \cite{FL16}*{Corollary~4.7}. 
Nevertheless, let us briefly explain why the argument remains unchanged in the twisted setting.

As functions of $g$, both sides of the above inequality are left $Q(\Q)N_{R^-}(\A)$-invariant. 
Hence, we may suppose $g=pak$, $p\in\omega\subset P_0(\A)\cap M_{R^-}(\A)$ (compact subset), $k\in K$, and $a\in A_0^G$ satisfies $F^Q(g,T) \, \tau_{P_0}^Q(H_{P_0}(a)-T_1) \, \tau_Q^{R^-}(H_Q(a)-T)=1$. 
There exists a compact set $D_{T_1}$ such that, for every $p\in\omega$, $k\in K$, and $a\in A_0^G$ satisfying $F^Q(g,T) \, \tau_{P_0}^Q(H_{P_0}(a)-T_1) \, \tau_Q^{R^-}(H_Q(a)-T)=1$, we have $p'\coloneqq a^{-1}pa\in D_{T_1}$.
Then,
\[
f(g^{-1}\eta\nu ng)=f(  y  \, (a^{-1}\nu n a) \, (p'k)), \qquad y\coloneqq k^{-1}a^{-1}p^{-1} \eta a.         
\]
Hence, by $p'k\in D_{T_1}K$ and reasoning as in the proof of \eqref{eq:2}, in order to prove the inequality, we can apply the proof of \cite{FL16}*{Proposition 4.6} to the function
\[
N_{S}^{M_{R^-}}(\A) \ni \nu \mapsto \int_{N_{R^-}(\A)} f(y n\nu)\, dn.
\]
Therefore, it suffices to index the coordinate families in \cite{FL16}*{Lemma~4.2} by $\Delta_{\tP_0}^{\tR^-}\setminus\Delta_{\tP_0}^{\tS}$.
\end{proof}

\subsection{Weighted parabolic integrals}
Take a Levi subset $\tM\in \cL^{\tG}(\tM_0)$ and a parabolic subset $\widetilde P\in\mathcal P^{\widetilde G}(\widetilde M)$. 
Set
\[
A_P^G\coloneqq A_M^G, \quad (A_P^G)^{1-\theta}\coloneqq \exp((\fa_P^G)^{1-\theta}),\quad  K_{M}\coloneqq M(\A)\cap K.
\]
\begin{proposition}\label{prop:biKP}
Let $\ell\ge0$ be an integer and let
$f\in\cC(\widetilde G(\A)^1)$ be bi-$K$-invariant. 
The integral
\begin{equation}\label{eq:biK-weighted-definition}
 f_{P,\ell}(x)
 =\int_{N_P(\A)}\int_{(A_P^G)^{1-\theta}}
       (1+\|H_P(a)\|)^\ell f(axn)\,da\,dn,\qquad x\in\widetilde M(\A)^1,
\end{equation}
converges absolutely at every $x$ and belongs to
$\cC(\widetilde M(\A)^1)$.
Note that $f_{P,\ell}$ is bi-$K_M$-invariant by the definition. 
For every $m\ge0$,
\begin{equation*}
 \|f_{P,\ell}\|_{\widetilde M(\A)^1,m}
 \le C_{P,\ell,m}\,
       \|f\|_{\widetilde G(\A)^1,m+\dim (\mathfrak a_P^G)^{\theta}}.
\end{equation*}
The constants depend only on the fixed data and the displayed indices,
not on $f$ or its support. 
\end{proposition}

We prove this proposition below.

\begin{lemma}\label{lem:biK-spherical-moment}
For every integer \(b\ge0\),
\begin{equation*}
\sup_{g\in G(\A)^1}
\int_K e^{\rho_P(H_P(kg))}
       (1+\|H_P(kg)\|)^b\,dk \ll_b 1.
\end{equation*}
\end{lemma}
\begin{proof}
If $P=G$, then $H_G(kg)=0$ for $g \in G(\A)^1$, and the assertion is immediate.
We therefore assume P is proper.
Normalize all compact Haar measures to have total volume one.

For any $\epsilon>0$, 
\[
(1+\|H\|)^b
\ll_{b,\epsilon}
\sum_{\alpha\in\Delta_P}\bigl(e^{\epsilon\alpha(H)}
           +e^{-\epsilon\alpha(H)}\bigr).
\]
Hence, it is sufficient to prove that there is a neighborhood $U$ of zero in
$(\fa_P^G)^*$ such that
\begin{equation}\label{eq:biK-exponential-average}
\sup_{g\in G(\A)^1}
\int_K e^{(\rho_P+\nu)(H_P(kg))}\,dk<\infty
\qquad(\nu\in U).
\end{equation}

For each place $v$, we normalize $\vol(K_v)=1$, and we can choose a minimal parabolic subgroup $P_{0,v}$ over $\Q_v$ contained in $P$.
Let $\rho_{0,v}$ and $\rho_{M,v}$ denote the half-sums of positive roots with multiplicities corresponding to $P_{0,v}$ and $P_{0,v}\cap M$.  
Let $\nu_v$ denote the pullback of $\nu\in(\fa_P^G)^*$ to $(\fa_{P_{0,v}}^G)^*$.
The local integral at $v$ in \eqref{eq:biK-exponential-average} is the spherical function
\[
\varphi_{\lambda,v}(g_v) \coloneqq  \int_{K_v}e^{(\lambda+\rho_{0,v})(H_{P_{0,v}}(k_vg_v))}\,dk_v, \quad g_v\in G(\Q_v) 
\]
with parameter $\lambda=-\rho_{M,v}+\nu_v$.

Consider the convex hull $\sC_v$ of the points $\{ w \rho_{0,v} \mid w\in W_{0,v}\}$ in $(\fa_{P_{0,v}}^G)^*$, where $W_{0,v}$ denotes the Weyl group of the root system at $v$.  
The support function $h_{\sC_v}(X)\coloneqq \max_{\lambda\in\sC_v}\lambda(X)$ satisfies 
\[
h_{\sC_v}(X)=\frac12\sum_{\alpha>0}\operatorname{mult}(\alpha)|\alpha(X)|, 
\]
where $\alpha$ runs over all positive roots and $\operatorname{mult}(\alpha)$ means the multiplicity of $\alpha$.
Then,
\[
h_{\sC_v}(X)-(-\rho_{M,v}(X))\ge \frac{1}{2}\sum_{\substack{\alpha>0\\ \alpha|_{\fa_M^G}\neq 0 }}|\alpha(X)|\ge 0 \quad \text{for all $X\in \fa_{P_{0,v}}^G$.}
\]
Since $\lambda\in \sC_v$ iff $\lambda(X)\le h_{\cC_v}(X)$ $(\forall X\in \fa_{P_{0,v}}^G)$, we have $-\rho_{M,v}\in\sC_v$.

We claim that a sufficiently small perturbation of
$-\rho_{M,v}$ in these directions still belongs to $\sC_v$.
To prove this, decompose the local root system into its
irreducible factors $\Sigma_{v,i}$.
Write $\Sigma_{M,v,i}$ for the Levi subsystem in $\Sigma_{v,i}$, and let $\sC_{v,i}$ denote the convex hull of the points $\{w_i\, \rho_{0,v,i} \mid w_i\in W_{0,v,i}\}$ in $(\fa_{P_{0,v}}^G)^*_i$, 
where the subscript $i$ denotes the corresponding component.
The convex hull $\sC_v$ is the product of the $\sC_{v,i}$ in the span of the local roots.

Suppose first that $\Sigma_{M,v,i}\subsetneq\Sigma_{v,i}$.
The roots outside the Levi subsystem span the whole root
space of this factor.
Indeed, their span is preserved by the simple reflections
belonging to the Levi subsystem, since these reflections
preserve the coefficients of the simple roots outside the Levi.
It is also preserved by each remaining simple reflection,
since its defining root belongs to this span.
The span is therefore a nonzero $W_{0,v,i}$-invariant
subspace, and irreducibility implies the assertion.
Consequently, on the corresponding coroot space,
\[
 q_{v,i}(X)
 \coloneqq
 \frac12
 \sum_{\alpha\in\Sigma_{v,i}^{+}
                    \setminus\Sigma_{M,v,i}^{+}}
       \operatorname{mult}(\alpha)|\alpha(X)|
\]
is positive for every $X\ne0$.
By compactness of the unit sphere, there is $c_{v,i}>0$ such that
\[
 h_{\sC_{v,i}}(X)-(-\rho_{M,v,i})(X)
 \ge q_{v,i}(X)
 \ge c_{v,i}\|X\|.
\]
Using the dual norm on the parameter space, we obtain
\[
 h_{\sC_{v,i}}(X)
       -(-\rho_{M,v,i}+\nu_{v,i})(X)
 \ge
 \bigl(c_{v,i}-\|\nu_{v,i}\|\bigr)\|X\|.
\]
Thus $-\rho_{M,v,i}+\nu_{v,i}\in\sC_{v,i}$ whenever $\|\nu_{v,i}\|<c_{v,i}$.

If instead $\Sigma_{M,v,i}=\Sigma_{v,i}$, then
$\nu_{v,i}=0$, because $\nu_v$ annihilates all the
coroots of this factor.
Moreover,
\[
 -\rho_{M,v,i}
 =-\rho_{0,v,i}
 =w_{0,v,i}\, \rho_{0,v,i}\in\sC_{v,i},
\]
where $w_{0,v,i}$ is the longest Weyl group element.
This proves that, for each $v$, there is a neighborhood $U_v$ of zero in $(\fa_P^G)^*$ such that
\[
 -\rho_{M,v}+\nu_v\in\sC_v
 \qquad(\nu\in U_v).
\]

The neighborhoods can be chosen uniformly in $v$.
Indeed, outside a finite set of places the local root data
with their Levi subsystems, root multiplicities, and
pullback maps $\nu\mapsto\nu_v$ have only finitely many
types: they arise from the fixed absolute root datum
by Frobenius actions through a finite group.
After making these identifications, only finitely many
of the preceding conditions on $\nu$ occur.
Intersecting the corresponding neighborhoods, and then
shrinking the intersection to accommodate the finitely
many remaining places, gives a single neighborhood $U$ of zero such that
\[
 -\rho_{M,v}+\nu_v\in\sC_v
 \qquad\text{for every place $v$ and every $\nu\in U$}.
\]

For classical boundedness criteria for spherical functions, see Helgason--Johnson \cite{HJ69} in the real case and Macdonald \cite{Mac71}*{\S4.7} in the non-archimedean case. 
The upper bound needed here follows directly from H\"older's inequality and Weyl invariance. 
Indeed, for $\nu\in U$ and each place $v$, write $\lambda_v=-\rho_{M,v}+\nu_v=\sum_{w\in W_{0,v}}t_w w\rho_{0,v}$, where $t_w\geq0$ and $\sum_w t_w=1$. 
Then
\[
0<\varphi_{\lambda_v,v}(g_v)
\leq
\prod_{\substack{w\in W_{0,v}\\
t_w>0}}\varphi_{w\rho_{0,v},v}(g_v)^{t_w}
=1
\qquad (g_v\in G(\mathbb{Q}_v)),
\]
since $\varphi_{-\rho_{0,v},v}=1$ and $-\rho_{0,v}\in W_{0,v}\rho_{0,v}$. Here the admissibility of $K$ ensures that each finite $K_v$ is a special maximal compact subgroup, so the Weyl invariance applies at every finite place.

For $g=(g_v)_v\in G(\A)^1$, the additivity of the height maps and the product decomposition of $K$ now give
\[
 \int_K e^{(\rho_P+\nu)(H_P(kg))}\,dk
 =
 \prod_v\varphi_{-\rho_{M,v}+\nu_v,v}(g_v)
 \le1
 \qquad(\nu\in U).
\]
For each fixed $g$, almost all $g_v$ belong to $K_v$;
the corresponding local integrands and spherical functions
are identically $1$.
Thus the displayed product has only finitely many
nontrivial factors.
This proves \eqref{eq:biK-exponential-average}, with upper bound $1$.
\end{proof}

\begin{proof}[Proof of Proposition~\ref{prop:biKP}]
Put $K'\coloneqq\theta(K)$. Since $f$ is bi-$K$-invariant, we have
\[
f_\delta(k_1gk_2')
=f(k_1g\delta\theta^{-1}(k_2'))
=f(g\delta)
=f_\delta(g)
\qquad (k_1\in K,\ k_2'\in K').
\]
Thus $f_\delta$ is left $K$-invariant and right $K'$-invariant.
For $X\in\cU(\fg_\C^1)$, $f_\delta*X=(f*\theta^{-1}(X))_\delta$. 
Since \(\theta^{-1}\) preserves the standard filtration
of the universal enveloping algebra, we have
\begin{equation}\label{eq:biK-pullback-seminorms}
\| f_\delta \|_{G(\A)^1,j}\asymp_j\|f\|_{\widetilde G(\A)^1,j}.
\end{equation}

Changing $n$ to $\theta(n)$,
we are reduced to studying
\begin{equation*}
F(y)\coloneqq
\int_{(\fa_P^G)^{1-\theta}}\int_{N_P(\A)}
(1+\|H\|)^\ell \, f_\delta(\exp(H)yn)\,dn\,dH,
\qquad y\in M(\A)^1,
\end{equation*}
with $f_{P,\ell}(y\delta)=F(y)$.
By \eqref{eq:biK-pullback-seminorms}, to prove the first assertion, it suffices to show that the integral $F(y)$ converges absolutely for every $y$ and that $F\in\cC(M(\A)^1)$.

For $X\in\mathcal U(\mathfrak g_\C^1)_{\le j}$, we set
\[
A_X(y)\coloneqq
\int_{(\fa_P^G)^{1-\theta}}\int_{N_P(\A)}
(1+\|H\|)^\ell \, \bigl|(f_\delta*X)(\exp(H)yn)\bigr|\,dn\,dH,
\qquad y\in M(\A)^1,
\]
\begin{lemma}\label{lem:biK-absolute-global}
\begin{equation*}
\int_{M(\A)^1}A_X(y)\,dy
\ll_{P,\ell,X} \|f_\delta\|_{G(\A)^1,j+\dim\fa_{\tP}^{\tG}}.
\end{equation*}    
\end{lemma}
\begin{proof}
Choose an $\Ad(K'_\infty)$-invariant Hermitian inner product on $\cU(\fg^1_\C)_{\le j+\dim\fa_{\tP}^{\tG}}$, and take an orthonormal basis $X_1,\dots,X_q$ of $\cU(\fg^1_\C)_{\le j+\dim\fa_{\tP}^{\tG}}$. 
Set $\cQ(g)\coloneqq (\sum_{j=1}^q|(f_\delta*X_j)(g)|^2)^{1/2}$. 
Fix a compact coordinate box $B\subset \fa_{\tP}^{\tG}$ containing zero in its interior.
Then,
\begin{itemize}
    \item $\cQ$ is left $K$-invariant and right $K'$-invariant.
    \item $\|\cQ\|_{L^1(G(\A)^1)}\ll_{P,j} \|f_\delta\|_{G(\A)^1,j+\dim\fa_{\tP}^{\tG}}$. Note that $q$ depends only on $j$ and $P$. 
    \item $(f_\delta*X)(\exp(H)yn)\ll_X \int_B\cQ(\exp(H)yn\exp(Y))\, dY $ is obtained by the local $L^1$ Sobolev inequality on $B$. 
\end{itemize}
Since $\rho_P\in (\fa_{\tP}^{\tG})^*$, we have $\rho_P(H^\theta)=\rho_P(H)$. 
Hence, $1\ll e^{\rho_P(H^{\theta})}=e^{\rho_P(H)}$ $(H^\theta\in B)$, and then
\begin{align*}
\int_{M(\A)^1}A_X(y)\,dy&\ll_{X} \int_{(\fa_P^G)^{1-\theta}\oplus B} \int_{M(\A)^1}\int_{N_P(\A)}(1+\|H^{1-\theta}\|)^\ell \cQ(\exp(H)yn)\, e^{\rho_P(H)} \, dn \, dy \, dH \\
&\ll_{X}  \int_{G(\A)^1}(1+\|H_P(g)\|)^\ell \cQ(g)\, e^{\rho_P(H_P(g))} \, dg\\
&\ll_{X} \int_{G(\A)^1}\cQ(g) \left(\int_K e^{\rho_P(H_P(kg))}(1+\|H_P(kg)\|)^\ell\,dk\right)dg.
\end{align*}
In passing to the integral over $G(\mathbb A)^1$, we use the Iwasawa decomposition with respect to $K'$, which is also admissible relative to $M_0$. 
Indeed, for $p\in P(\mathbb A)\cap G(\mathbb A)^1$ and $k'\in K'$, we have $Q(pk')=Q(p)$ and $H_P(pk')=H_P(p)+H_P(k')$. 
Since $H_P(K')$ is bounded, the weight
\[
w(H)=(1+\|H\|)^\ell e^{\rho_P(H)}
\]
satisfies \(w(H_P(pk'))\asymp w(H_P(p))\) uniformly in \(p\) and \(k'\). The last line then follows by averaging on the left over $K$, using the left $K$-invariance of $Q$.
By Lemma~\ref{lem:biK-spherical-moment} and $\|\cQ\|_{L^1(G(\A)^1)}\ll_{j} \|f_\delta\|_{G(\A)^1,j+\dim\fa_{\tP}^{\tG}}$, 
\[
\int_{M(\A)^1}A_X(y)\,dy 
\ll_{\ell}\|\cQ\|_{L^1(G(\A)^1)} \ll_{P,j} \|f_\delta\|_{G(\A)^1,j+\dim\fa_{\tP}^{\tG}}.
\]
This completes the proof.
\end{proof}
\begin{lemma}\label{lem:biK-absolute-pointwise}
\begin{equation*}
\sup_{y\in M(\A)^1}A_X(y)
\ll_{P,\ell,X}\|f_\delta\|_{G(\A)^1,j+\dim M(\R)^1 + \dim\fa_{\tP}^{\tG}}.
\end{equation*}
\end{lemma}
\begin{proof}
Fix an identity neighborhood $U_M\subset M(\R)^1$.
Applying the local $L^1$ Sobolev inequality to $U_M\ni u\mapsto(f_\delta*X)(\exp(H)ynu)$, 
\[
A_X(y) \ll_X \sum_{X'}\int_{U_M}A_{X'}(yu)\,du,
\]
where the sum is finite $(\le \dim \cU(\fg^1_\C)_{\le j+\dim M(\R)^1})$ and $\deg X'\le j+\dim M(\R)^1$.
Since the archimedean differential operators preserve right $K'_f$-invariance and conjugation by $K'_f\cap M(\mathbb A_f)$ preserves Haar measure on $N_P(\mathbb A)$, each $A_{X'}$ is right invariant under $K'_f\cap M(\mathbb A_f)$. 
The assertion now follows from Lemma~\ref{lem:biK-absolute-global}.
\end{proof}

By Lemmas \ref{lem:biK-absolute-global} and \ref{lem:biK-absolute-pointwise},  the integral $F(y)$ converges absolutely at every $y\in M(\A)^1$, and belongs to $L^1(M(\A)^1)$.
Put $I_X(y)\coloneqq
\int_{(\fa_P^G)^{1-\theta}}\int_{N_P(\A)}
(1+\|H\|)^\ell \, (f_\delta*X)(\exp(H)yn)\,dn\,dH$. 
Set $\fm^1\coloneqq \mathrm{Lie}(M(\R)^1)$. 
For any $X\in\cU(\fm^1_\C)\subset \cU(\fg^1_\C)$, we have 
\[
 |I_X(y)|\le A_X(y).
\]
Hence, by Lemma \ref{lem:biK-absolute-pointwise}, the integral defining $I_X(y)$ converges absolutely at every $y\in M(\A)^1$.
After a measure-preserving change of variables in $N_P(\A)$, local Sobolev estimates and Lemma~\ref{lem:biK-absolute-pointwise} provide locally uniform integrable majorants for the right derivatives of the integrand.
Thus, differentiation under the integral sign is justified by dominated convergence, giving
\[
(F*X)(y)=I_X(y).
\]
Therefore, we obtain $F\in\cC(M(\A)^1)$ by Lemmas \ref{lem:biK-absolute-global} and \ref{lem:biK-absolute-pointwise}. 
In addition, the second assertion also follows from Lemma \ref{lem:biK-absolute-global}. 
\end{proof}

\section{On the decay of the individual terms for fixed compact support}
\label{sec:cga}

We retain the notation and measures of \S~\ref{sec:preliminaries}.
Fix a compact subset $\Omega\subset\tG(\A)^1$ and an open compact
subgroup $K_0\subset G(\A_f)$. 
In this section, we will study the decay properties for test functions in $C_\Omega^\infty(\tG(\A)^1;K_0)$, which was defined in \S~\ref{sec:Frechet}.

\subsection{The main term}\label{sec:maincoarse}

Let $Q$, $R\in\cF^G(M_0)$ such that $P_0\subset Q\subset R$. 
For $B_{\tfo,Q,R}^T(x,f)$, the alternating sum has one term precisely when $Q^+=R^-$.
We retain all such terms in the main kernel:
\begin{equation*}
 m_{\tfo}^T(x,f)
   \coloneqq \sum_{\substack{Q\subset R\\Q^+=R^-}}
         B_{\tfo,Q,R}^T(x,f).
\end{equation*}
If $Q\neq G$, then $\widetilde\sigma_Q^Q=0$, and so $B_{\tfo,Q,Q}^T(x,f)=0$. 
Hence, 
\begin{align*}    
 m_{\tfo}^T(x,f)=& F^G(x,T)K_{\tG,\tfo}(x,x) \\
 &+\sum_{\substack{Q\subsetneq R\\Q^+=R^-=P}}
          (-1)^{\dim\fa_{\tP}^{\tG}}\sum_{\xi\in Q(\Q)\backslash G(\Q)}F^Q(\xi x,T) \widetilde\sigma_Q^R(H_0(\xi x)-T)
                 K_{\tP,\tfo}(\xi x,\xi x).
\end{align*}
Set
\[
 M_{\tfo}^T(f)
  \coloneqq\int_{G(\Q)\bs G(\A)^1}m_{\tfo}^T(x,f)\,dx,
\]
and then we have
\[
 M_{\tfo}^T(f)= \int_{G(\Q)\bs G(\A)^1}F^G(x,T)K_{\tG,\tfo}(x,x)\,dx+\sum_{\substack{Q\subsetneq R\\Q^+=R^-=P}}
          (-1)^{\dim\fa_{\tP}^{\tG}}\, M_{\tfo,Q,R}^T(f),
\]
where
\begin{equation*}
 M_{\tfo,Q,R}^T(f)
 \coloneqq \int_{Q(\Q)\bs G(\A)^1}F^Q(x,T)\widetilde\sigma_Q^R(H_0(x)-T)
                 K_{\tP,\tfo}(x,x)\,dx.
\end{equation*}
The sign $(-1)^{\dim\fa_{\tP}^{\tG}}$ is outside this integral.

\subsection{Some lemmas}

\begin{itemize}
    \item We set $W^M_0 \coloneqq N_{M}(T_0)(\Q)/M_0(\Q)$ for $M\in\cL^G(M_0)$. 
    \item The projection onto $(\fa_{M_0}^G)^\theta$ is $\pi_\theta\coloneqq \frac1\ell\sum_{j=0}^{\ell-1}\theta^j$. 
          For $H\in\fa_{M_0}^G$, we set $H^\theta\coloneqq \pi_\theta(H)$ and denote by $H^-$ its complementary component. 
          Then, 
          \[
          H=H^\theta+H^-, \quad  H^\theta\in\fa_{\tM_0}^{\tG}, \quad H^-\in (\fa_{M_0}^G)^{1-\theta}.
          \]
          The superscript $-$ does not denote the $(-1)$-eigenspace.
\end{itemize}

\begin{lemma}\label{lem:fsa-compact-support}
Let $\tP\in\cF^{\tG}(\tM_0)$ such that $P_0\subset P$.
Put $M=M_P$ and $N_{P_0}^M=N_{P_0}\cap M$.
Fix a compact subset $C_M\subset M(\A)$ and $c\in\R$.
If $a=\exp H\in A_0^G$, $\alpha(H)\geq c$ for all $\alpha\in \Delta_{P_0}^P$, and $\eta\in M(\Q)$, $n_1,n_2\in N_{P_0}^M(\A)$ satisfy $a^{-1}n_1\eta n_2\theta(a)\in C_M$, then 
\[
\|(1-\theta)H\|+\|H^-\|\ll_{C_M,c}1.
\]
\end{lemma}
\begin{proof}
Suppose $a^{-1}n_1\eta n_2\theta(a)\in C_M$. 
Fix a representative $\dot w\in N_M(M_0)(\Q)$ for $w\in W_0^M$.
The Bruhat decomposition in $M$ gives
\[
  \eta=u_\eta m_\eta\dot w v_\eta,
  \qquad
  u_\eta,v_\eta\in N_{P_0}^M(\Q),
  \quad m_\eta\in M_0(\Q),
\]
for some $w\in W_0^M$.
Put $H\coloneqq H_0(a)$ and $v'\coloneqq \theta(a)^{-1}v_\eta n_2\theta(a)$.
Then,
\[
  H_0(a^{-1}n_1\eta n_2\theta(a))=-H+w\theta H+H_0(\dot w v').
\]
By \cite{LW13}*{Lemma~3.3.1(i)}, there exists a constant $B>0$ such that  
\[
  \varpi\bigl(H_0(\dot w v)\bigr)\le B, \qquad w\in W_0^M
  \qquad v\in N_{P_0}^M(\A),
  \quad \varpi\in\widehat\Delta_{P_0}^P.
\]
Note that $B$ is independent of $v$. 
By the compactness of $C_M$, there exists a constant $B(C_M)>0$ such that $|\varpi\bigl(H_0(a^{-1}n_1\eta n_2\theta(a))\bigr)|\le B(C_M)$. 
Hence, for any $\varpi\in\widehat\Delta_{P_0}^P$,
\begin{equation}\label{eq:wH1}
\varpi(H)=-\varpi(H_0(a^{-1}n_1\eta n_2\theta(a)))+\varpi(w\theta H)+\varpi(H_0(\dot w v'))\le B(C_M)+B+\varpi(w\theta H).    
\end{equation}

Suppose that $\alpha(H)\geq c$ holds for every $\alpha\in \Delta_{P_0}^P$. 
Let $H_c\in\fa_{P_0}^P$ be the unique vector satisfying $\alpha(H_c)=c$, $\forall\alpha\in\Delta_{P_0}^P$. 
The automorphism $\theta$ permutes $\Delta_{P_0}^P$, so
$\theta H_c=H_c$. 
Let $H^P$ denote the projection of $H$ into $\fa_{P_0}^P$. 
Then, $X\coloneqq \theta H^P-H_c\in\fa_{P_0}^P$ satisfies $\alpha(X)\ge 0$, $\forall\alpha\in\Delta_{P_0}^P$.
Hence, the usual dominance property of Weyl groups gives
\[
  X-wX\in
  \sum_{\alpha\in\Delta_{P_0}^P}\R_{\geq 0}\alpha^\vee.
\]
Every $\varpi\in\widehat\Delta_{P_0}^P$ is nonnegative on this
cone. Since $W_0^M$ acts trivially on $\fa_P^G$, it follows that
\begin{equation}\label{eq:wH2}  \varpi(w\theta H-\theta H)
    =\varpi(wX-X)+\varpi(wH_c-H_c)
    \leq\varpi(wH_c-H_c)\leq B_c
\end{equation}
for a constant $B_c\geq0$ depending only on $c$ and the fixed root
data. 
By \eqref{eq:wH1} and \eqref{eq:wH2},
\[
  \varpi\bigl((1-\theta)H\bigr)\leq B',
  \qquad \varpi\in\widehat\Delta_{P_0}^P,
  \qquad B'\coloneqq B+B(C_M)+B_c.
\]

Since $\tvarpi\in \widehat\Delta_{\tP_0}^{\tP}$ satisfies $\tvarpi\bigl((1-\theta)H\bigr)=0$ and $\tvarpi=\frac{1}{m+1}(\varpi+\sum_{j=1}^m \varpi_j)$ $(\exists\varpi_1,\dots,\varpi_m\in \widehat\Delta_{P_0}^P)$, we obtain 
\[
\varpi\bigl((1-\theta)H\bigr)\ge -\sum_{j=1}^m \varpi_j\bigl((1-\theta)H\bigr)\ge -m B'.
\]
As $\widehat\Delta_{P_0}^P$ is a basis of
$(\fa_{P_0}^P)^*$, these two-sided bounds imply $\bigl\|((1-\theta)H)^P\bigr\|\ll_{C_M,c}1$.

It remains to bound the $\fa_P^G$-component. 
Since $H_P$ vanishes on $M(\A)^1$, 
\[
  H_P(a^{-1}n_1\eta n_2\theta(a))=-H_P(a)+\theta H_P(a)=-(1-\theta)H_P(a).
\]
Since $a^{-1}n_1\eta n_2\theta(a)\in C_M$, this bounds $(1-\theta)H_P$.
Thus, we conclude that $\|(1-\theta)H\|\ll_{C_M,c}1$.

Recall that $H^-$ is the component orthogonal
to the $\theta$-fixed subspace. Thus
\[
  H^-=H-\frac1\ell\sum_{j=0}^{\ell-1}\theta^jH
      =\frac1\ell\sum_{j=1}^{\ell-1}
         \sum_{k=0}^{j-1}\theta^k(1-\theta)H.
\]
This explicitly expresses the inverse of $1-\theta$ on the
complement of its kernel and gives
$\|H^-\|\ll_\theta\|(1-\theta)H\|$.
Hence, the proof is completed. 
\end{proof}

\begin{lemma}\label{lem:cga-lattice}
Let $\tP\in\cF^{\tG}(\tM_0)$ containing $\tP_0$, $M\coloneqq M_P$, and $\tM(\A)^{(1)}\coloneqq \tM(\A)\cap\tG(\A)^1$. 
Fix a compact subset $\Omega_M=\Omega_{M,\infty}\times \Omega_{M,0}\subset M(\A)$, where $\Omega_{M,\infty}\subset M(\R)$ and $\Omega_{M,0}\subset M(\A_f)$.
We also fix a constant $c\in\R$ and a continuous function $f$ on $\tM(\A)^{(1)}$ whose support is contained in $(\Omega_M \delta)\cap \tM(\A)^{(1)}$. 
If $a=\exp H\in A_0^G$, $\alpha(H)\geq c$ for all $\alpha\in \Delta_{P_0}^P$, then one has
\[
\sum_{\eta\in \tM(\Q)} |f(  a^{-1}\eta a)|\ll_{\Omega_M,c} e^{2\rho_{P_0\cap M}(H)} \|f\|_\inf    .
\]
\end{lemma}
\begin{proof}
We first work on $M(\A)$ by replacing $f$ by $f_\delta$. 
Suppose that $a=\exp H\in A_0^G$ and $\alpha(H)\geq c$ $(\forall \alpha\in \Delta_{P_0}^P)$. 
We first establish a counting estimate with fixed finite support.
Set $J\coloneqq K(1)\cap M(\A_f)$, see \S\ref{sec:level} for $K(1)$.
For any element $x\in \Omega_{M,0}$, we claim that
\begin{equation}\label{eq:cga-lattice-fixed-support}
\mathcal{N}_x(a)\coloneqq \#\{\gamma\in M(\Q)\mid
       a^{-1}\gamma \theta(a)\in \Omega_{M,\inf}\times xJ\}
 \ll_{\Omega_M,c}e^{2\rho_{P_0\cap M}(H)}.
\end{equation}
By Lemma \ref{lem:fsa-compact-support}, we may absorb the bounded component $H^-$ into a fixed enlargement of $\Omega_{M,\infty}$ and hence assume $H=H^\theta$. This replacement preserves the root inequalities and the value of $\rho_{P_0\cap M}(H)$.
Hence, 
\[
\mathcal{N}_x(a)\ll_{\Omega_M,c} \#\{\gamma\in M(\Q)\mid
       a^{-1}\gamma a\in \Omega_{M,\inf}\times xJ\}.
\]
Set
\[
 D_c\coloneqq
 \{H\in\fa_{M_0}^G\mid
     \alpha(H)\geq c\quad(\alpha\in\Delta_{P_0}^P)\}.
\]
For sufficiently large $b\, (>c)$ and $P_0\subseteq Q\subseteq P$, put
\[
 D_{b,Q}\coloneqq
 \left\{H\in D_c\; \middle| \;
 \begin{array}{ll}
 \alpha(H)\leq b
   &(\alpha\in\Delta_{P_0}^Q),\\
 \alpha(H)>b
   &(\alpha\in\Delta_{P_0}^P\setminus\Delta_{P_0}^Q)
 \end{array}\right\}.
\]
The correspondence between standard parabolic subgroups and
subsets of $\Delta_{P_0}$ gives
\[
 D_c=\bigsqcup_{P_0\subseteq Q\subseteq P}D_{b,Q}.
\]
Since $b$ is sufficiently large, we have
\[
 H\in D_{b,Q},\quad
 a^{-1}\gamma a\in \Omega_{M,\infty}\times xJ
 \quad\Longrightarrow\quad
 \gamma\in(Q\cap M)(\Q).
\]
Note that $x$ belongs to $\Omega_{M,0}$, so $b$ depends on $\Omega_M$, $J$, and $c$. 
Assume that $a^{-1}\gamma a\in \Omega_{M,\infty}\times xJ$.
Fix $Q$, and suppose $H\in D_{b,Q}$. 
Then, we can write
\[
 \gamma=\mu u,\qquad
 \mu\in M_Q(\Q),\quad u\in(N_Q\cap M)(\Q).
\]
Then, $a^{-1}\mu a$ and $a^{-1}ua$ belong to fixed compact
sets in $M_Q(\A)$ and $(N_Q\cap M)(\A)$, respectively.
Since the projection of $H$ onto $\fa_{M_0}^Q$ is bounded,
$\mu$ belongs to a fixed compact subset of $M_Q(\A)$, which contains only finitely many elements of
$M_Q(\Q)$.
On the other hand, the number of possible $u$ is bounded by a constant multiple of $e^{2\rho_{Q\cap M}(H)}$. 
Therefore,
\[
       \#\{\gamma\in M(\Q)\mid
       a^{-1}\gamma \theta(a)\in \Omega_{M,\infty}\times xJ\} \ll_{\Omega_M,c} \sum_{P_0\subseteq Q\subseteq P}e^{2\rho_{Q\cap M}(H)} \ll e^{2\rho_{P_0\cap M}(H)}.
\]
Since the implicit constant does not depend on $x$, the proof of \eqref{eq:cga-lattice-fixed-support} is completed.

There are $x_1,\dots,x_q\in \Omega_{M,0}$ such that $\Omega_{M,0}$ is covered by $ \cup_{j=1}^q  x_j J$. 
Hence, by \eqref{eq:cga-lattice-fixed-support} we obtain
\[
\sum_{\eta\in\tM(\Q)}|f(a^{-1}\eta a)| \le   \|f \|_\inf \sum_{ j=1 }^q  \, \mathcal{N}_{x_j}(a)\ll_{\Omega_M,c} e^{2\rho_{P_0\cap M}(H)}\, \|f\|_\inf  \qedhere
\]
\end{proof}

\begin{lemma}\label{lem:cga-cone-fibers}
Let $Q$, $R\in\cF^G(M_0)$, suppose $P_0\subset Q\subset R$, $Q^+\subset R^-$. 
There is a constant $c_{Q,R}$ such that
\begin{equation*}
 H\in \fa_{\tR^-}^{\tG},\quad \widetilde\sigma_Q^R(H+Z)=1
       \quad\Longrightarrow\quad \|H\|\leq c_{Q,R}\, \|Z\|
\end{equation*}
for every $Z\in \fa_{M_0}^G$.
\end{lemma}
\begin{proof}
See \cite{LW13}*{Lemma 2.11.6 (ii)}.
\end{proof}

Set 
\[
\beta_{Q,R}(H)\coloneqq \sum_{\tbeta\in
       \Delta_{\tP_0}^{\tR^-}\setminus\Delta_{\tP_0}^{\tQ^+}}
       \ \min_{\substack{\alpha\in\Delta_{P_0}\\
                         \talpha=\tbeta}}
       \alpha_{Q}(H) ,\qquad  H\in \fa_{M_0}^G.
\]
For $c\geq0$, define
\begin{align*}
 \mathcal D_c^{Q,R}(T)\coloneqq \{H\in \fa_{M_0}^G \mid & \ \|H^-\|\leq c,\quad
       \alpha(H-T_1)>0\;\; (\forall\alpha\in \Delta_{P_0}^Q), \\
  &    \ \varpi(H-T)\leq0\;\; (\forall\varpi\in\widehat\Delta_{P_0}^Q),\quad \widetilde\sigma_Q^R(H-T)=1\}.    
\end{align*}
\begin{lemma}\label{lem:fsa-cone-integral}
Let $Q$, $R\in\cF^G(M_0)$, suppose $P_0\subset Q\subset R$, $Q^+\subset R^-$. 
For $c\ge 0$, $\kappa\geq0$ and $T\in \fa_{M_0}^G$,
\begin{equation*}
 \int_{\cD_c^{Q,R}(T)}
       (1+\|H\|)^\kappa e^{-\beta_{Q,R}(H-T)}\,dH
 \ll_{c,\kappa}
       (1+\|T\|)^{\dim \fa_{\tM_0}^{\tG}-\dim\fa_{\tQ^+}^{\tR^-}+\kappa}.
\end{equation*}
\end{lemma}
\begin{proof}
Suppose $H\in \cD_c^{Q,R}(T)$. 
For $H=H^\theta+ H^-$, we have $\|H^-\|\le c$.  
Following $\fa_{\tM_0}^{\tG}=\fa_{\tP_0}^{\tQ^+}\oplus\fa_{\tQ^+}^{\tR^-}\oplus\fa_{\tR^-}^{\tG}$, we write 
\[
H^{\theta}=H_1+H_2+H_3, \quad H_1\in \fa_{\tP_0}^{\tQ^+}, \quad H_2\in\fa_{\tQ^+}^{\tR^-}, \quad H_3\in\fa_{\tR^-}^{\tG}.
\]
By the conditions $\alpha(H-T_1)>0$ $(\forall\alpha\in \Delta_{P_0}^Q)$ and $\varpi(H-T)\leq0$ $(\forall\varpi\in\widehat\Delta_{P_0}^Q)$, 
\[
\|H_1\| \ll_c 1+\|T\|.
\]
In addition, by Lemma \ref{lem:cga-cone-fibers}, 
\[
\|H_3\|\ll_c  1+\|H_2\|+\|T\|.
\]
Hence, 
\begin{equation}    \label{eq:1019e1}
 \int_{\cD_c^{Q,R}(T)}
       (1+\|H\|)^\kappa e^{-\beta_{Q,R}(H-T)}\,dH  \ll_{c,\kappa} \int_{\cD_c^{Q,R}(T)} e^{-\beta_{Q,R}(H-T)}
       (1+\|H_2\|+\|T\|)^{\kappa} \, dH.    
\end{equation}

For $\tbeta\in \Delta_{\tP_0}^{\tR^-}\setminus\Delta_{\tP_0}^{\tQ^+}$, we set $t_{\tbeta}\coloneqq \min_{\substack{\alpha\in\Delta_{P_0}\\ \talpha=\tbeta}} \alpha_{Q}(H-T)$, then
\[
\beta_{Q,R}(H-T)=\sum_{\tbeta\in
       \Delta_{\tP_0}^{\tR^-}\setminus\Delta_{\tP_0}^{\tQ^+}}t_{\tbeta} , \qquad t_{\tbeta}=\tbeta_Q(H_2)+\min_{\substack{\alpha\in\Delta_{P_0}\\ \talpha=\tbeta}}\alpha_Q(H_1+H^- -T).
\]
By $\widetilde\sigma_Q^R(H-T)=1$, we have $t_{\tbeta}>0$ for all $\tbeta\in \Delta_{\tP_0}^{\tR^-}\setminus\Delta_{\tP_0}^{\tQ^+}$. 
In addition,
\[
\|H_2\|\ll \sum_{\tbeta\in\Delta_{\tP_0}^{\tR^-}\setminus\Delta_{\tP_0}^{\tQ^+}}|\tbeta_Q(H_2)|\ll_c 1+\|T\|+\sum_{\tbeta\in\Delta_{\tP_0}^{\tR^-}\setminus\Delta_{\tP_0}^{\tQ^+}} t_{\tbeta}.
\]
Thus, on the right-hand side of \eqref{eq:1019e1}, integrating first with respect to $H_3$, then with respect to $H_2$ (the change of variables to $t_{\tbeta}$), and finally with respect to $H_1$ and $H^-$, we obtain
\begin{align*}
& \int_{\cD_c^{Q,R}(T)}
       (1+\|H_2\|+\|T\|)^{\kappa} e^{-\beta_{Q,R}(H-T)}\,dH  \\
&\ll_{c} \int_{\substack{t_{\tbeta}>0 \;\; (\forall \tbeta)\\ \|H_1\|\ll 1+\|T\|\\ \|H^-\|\le c}}
       (1+\|H_2\|+\|T\|)^{\kappa+\dim\fa_{\tR^-}^{\tG}} e^{-\beta_{Q,R}(H-T)}\,dH_2\, dH_1\, dH^- \\
& \ll_\kappa (1+\|T\|)^{\kappa+\dim\fa_{\tR^-}^{\tG}}  \int_{\substack{\|H_1\|\ll 1+\|T\|\\ \|H^-\|\le c}} \, dH_1\, dH^-\ll (1+\|T\|)^{\dim \fa_{\tM_0}^{\tG}-\dim\fa_{\tQ^+}^{\tR^-}+\kappa}.        
\end{align*}
This completes the proof.
\end{proof}

\subsection{Approximation theorem}\label{sec:apprpx}

\begin{theorem}\label{thm:cga}
There exist integers $m$, $s$, $r\geq0$ such that for each compact $\Omega\subset\tG(\A)^1$ there
are constants $C_\Omega$, $t_\Omega>0$ with the following property:
for every open compact subgroup $K_0\subset G(\A_f)$,
\begin{align*}
& \sum_{\tfo\in\tcO}
    |J_{\tfo}^{T,\tG}(f)-M_{\tfo}^T(f)|
 \leq\sum_{\tfo\in\tcO}\int_{G(\Q)\bs G(\A)^1}
          |k_{\tfo}^T(x,f)-m_{\tfo}^T(x,f)|\,dx\\
&\leq C_\Omega\level(K_0)^s
           q_m(f)(1+\|T\|)^r e^{-d(T)}    
\end{align*}
for $f\in\cC_\Omega^\infty(\tG(\A)^1;K_0)$ and $d(T)\geq t_\Omega$.
The constants $C_\Omega$, $t_\Omega$ are independent of $K_0$, $f$,
and $T$; the integers $m$, $s$, $r$ are also independent of $\Omega$.
\end{theorem}
\begin{remark}\label{rem:untwistedmain}
In the untwisted case $\tG=G$ and $\delta=1$, the only nonzero one-term contribution is
$(Q,R)=(G,G)$, so this theorem reduces to the result \cite{FL16}*{Theorem 5.1 (3)}.
\end{remark}
\begin{proof}
By unfolding the right hand side of
\[
\int_{G(\Q)\bs G(\A)^1} | k_{\tfo}^T(x,f)-m_{\tfo}^T(x,f)| \, dx
   \le \sum_{\substack{Q\subset R\\Q^+ \subsetneq R^-}}
         \int_{G(\Q)\bs G(\A)^1} |B_{\tfo,Q,R}^T(x,f)|\, dx
\]
it is enough to bound
\[
 \cR_{Q,R}(T,f)
  \coloneqq\sum_{\tfo\in\tcO}\int_{Q(\Q)\bs G(\A)^1}
      F^Q(x,T)\widetilde\sigma_Q^R(H_0(x)-T)
      \Bigl|\sum_{\tQ^+\subset \tP\subset \tR^-} (-1)^{\dim\fa_{\tP}^{\tG}}  K_{\tP,\tfo}( x, x)\Bigr|\,dx
\]
for each pair $(Q,R)$ such that $\tQ^+\subsetneq \tR^-$.
Recall that $d(T)>t_\Omega$, with $t_\Omega$ sufficiently large, and $F^Q(x,T)\, \tsigma_Q^R(H_0(x)-T)=1$ implies $\tau_Q^{R^-}(H_0(x)-T)=1$. 
Hence, by the compactness of the support of $f$, 
\[
\sum_{\tQ^+\subset \tP\subset \tR^-} (-1)^{\dim\fa_{\tP}^{\tG}}  K_{\tP,\tfo}( x, x)=\cE_{Q,R;Q^+}(x,\tM_{Q^+}(\Q)\cap\tfo,f).
\]
Applying Proposition \ref{prop:tfl47} to $\cE_{Q,R;Q^+}(x,\tM_{Q^+}(\Q)\cap\tfo,f)$, we obtain
\begin{align*}
\cR_{Q,R}(T,f)\ll& \level(K_0)^s\sum_{\tfo\in\tcO}\sum_{\gamma\in \tM_{Q^+}(\Q)\cap\tfo}  
      \int_{Q(\Q)\bs Q(\A)\cap G(\A)^1}\int_K F^Q(q,T)\widetilde\sigma_Q^R(H_0(q)-T) \\
      &  e^{-\Gamma_{Q,R;Q^+}(T,H_0(q))} \sum_{j=1}^m\int_{N_{Q^+}(\A)} |f_j(k^{-1}q^{-1}\gamma nqk)|\,dn\, dk \, dq ,
\end{align*}
where $f_j\coloneqq f*X_j$. 
Since it suffices to obtain an upper bound for the integral, we enlarge the domain of integration in $q$ to the set $\cS^Q(T_1,\omega)$, see Remark \ref{rem:truncation-Arthur-compatibility}. 
Set $M\coloneqq M_{Q^+}$ and recall $|f_{j}|_{Q^+}(x)\coloneqq \int_{N_{Q^+}(\A)} \int_K |f_j(k^{-1}xn k)| \, dk \, dn$, $x\in \tM(\A)^{(1)}$. 
Then, 
\begin{align*}
\cR_{Q,R}(T,f)\ll& \level(K_0)^s \sum_{j=1}^m  
      \int_{A_0^G} \int_{\omega\cap N_{P_0}^M(\A)} F^Q(a,T)\, \widetilde\sigma_Q^R(H_0(a)-T) \\
      &  \sum_{\gamma\in M(\Q)}e^{-\Gamma_{Q,R;Q^+}(T,H_0(a))-2\rho_{P_0}(H_0(a))+ 2\rho_{Q^+}(\theta H_0(a))}  |f_{j}|_{Q^+}(a^{-1}u^{-1}\gamma\theta(ua) \delta) \, du \, da.
\end{align*}
Put $H\coloneqq H_0(a)$.  
Suppose $F^Q(a,T)\, \widetilde\sigma_Q^R(H_0(a)-T)=1$ and $|f_{j}|_{Q^+}(a^{-1}u^{-1}\gamma\theta(ua)\delta)\neq 0$. 
Then, by Lemma \ref{lem:fsa-compact-support}, we have 
\[
\|(1-\theta)H\|+ \|H^-\|<c_\Omega
\]
for some constant $c_\Omega>0$. 
In addition, 
\[
-\Gamma_{Q,R;Q^+}(T,H)-2\rho_{P_0}(H)+ 2\rho_{Q^+}(\theta H)\le  -d(T)-\beta_{Q,R}(H-T)-2\rho_{P_0\cap M}(H).
\]
Therefore, by applying Lemma \ref{lem:cga-lattice} to it,
\[
\cR_{Q,R}(T,f)\ll_\Omega \level(K_0)^{s} e^{-d(T)}  \sum_{j=1}^m  \||f_j|_{Q^+}\|_{\infty}  \int_{\cD^{Q,R}_{c_\Omega}(T)}    e^{-\beta_{Q,R}(H-T)}  \, dH,
\]
since the sum for $\tfo$ is finite and the range of $a^{-1}ua$ is compact. 
Thus, the proof is completed by Lemma \ref{lem:fsa-cone-integral} and $\||f_j|_{Q^+}\|_{\infty}\ll_\Omega \|f_j\|_\inf$. 
\end{proof}

\subsection{The one-term contributions}

Assume that we are in one of the following two cases.
\begin{itemize}
    \item $Q=P=R=G$.
    \item $Q\subsetneq R$ and $Q^+=P=R^-$. (This case contains the case $R=G$.)
\end{itemize}
\begin{proposition}\label{prop:onetermconvergence}
There exist integers $m$, $r$, $s\geq0$ such that for each compact $\Omega\subset\tG(\A)^1$ there
are constants $C_\Omega$, $t_\Omega>0$ with the following property:
for every open compact subgroup $K_0\subset G(\A_f)$,
\[
\int_{G(\Q)\bs G(\A)^1} \Bigl| B^T_{\fo,Q,R}(x,f) \Bigr|\, dx \ll_{\Omega}\level(K_0)^r\,  (1+\|T\|)^s \, q_m(f) .
\]
for $f\in\cC_\Omega^\infty(\tG(\A)^1;K_0)$ and $d(T)\geq t_\Omega$.
\end{proposition}
\begin{proof}
This follows from the same argument as in the proof of Theorem \ref{thm:cga}.
\end{proof}
We set $M\coloneqq M_P$, $N_Q^M\coloneqq N_Q\cap M$, $A_Q^G\coloneqq \{ a\in A_{Q} \mid H_G(a)=0 \}$, and $f_{P}(x)\coloneqq \int_K \int_{N_P(\A)}f(k^{-1} x u k)\,du \, dk$, $x\in\tM(\A)^{(1)}=\tM(\A)\cap\tG(\A)^1$. 
It follows from the proof of Proposition \ref{prop:onetermconvergence} that
\begin{align*}
     M_{\tfo,Q,R}^T(f)
 ={}&\int_{A_Q^G}\int_{M_Q(\Q)\backslash M_Q(\A)^1} \int_{N_Q^M(\Q)\backslash N_Q^M(\A)} F^Q(m,T)\, \tsigma_Q^R(H_0(a)-T)      e^{-2\rho_{Q\cap M}(H_0(a))}\\
    &    \sum_{\gamma\in\tM(\Q)\cap \tfo} f_{P}(a^{-1}m^{-1}n^{-1}\gamma nma )
          \,dn \, dm \, da.
\end{align*}
This formula gives a parabolic descent to $\tM$ with test function $f_P$.
Since the $\theta$-saturation of $Q\cap M$ in $M$ is $M$ and $R\cap M=M$, the inner kernel has the same form as in the case $R=G$, with $G$ replaced by $M$.
The integration over $A_M^G$ and the factor $\widetilde\sigma_Q^R$ remain to be taken into account.

To see that these terms do not in fact decay as $d(T)$ increases, 
it is sufficient to consider the case $R=G$, hence $Q\subsetneq Q^+=P=R^-=G$. 
Set $f_K(x)\coloneqq \int_{K}f(k^{-1}xk)d k$, $x\in\tG(\A)^1$. 
Then, 
\begin{align}\label{eq:R=G}
     M_{\tfo,Q,G}^T(f)
 ={}&\int_{A_Q^G}\int_{M_Q(\Q)\backslash M_Q(\A)^1} \int_{N_Q(\Q)\backslash N_Q(\A)} F^Q(m,T)\, \tau_Q(H_0(a)-T)
                           e^{-2\rho_{Q}(H_0(a))}\\
    &
        \sum_{\gamma\in\tfo} f_K(a^{-1}m^{-1}n^{-1}\gamma nma)
          \,dn \, dm \, da. \nonumber \\
    ={}& \int_{G(\Q)\bs G(\A)^1} K_{\tG,\tfo}(x,x) \sum_{\xi\in Q(\Q)\bs G(\Q)} F^Q(\xi x,T)\tau_Q(H_0(\xi x)-T)\, dx. \label{eq:R=G K_G}
\end{align}
Let $a=\exp(H)\in A_0^G$, and $x=pak\in \cS^G(T_1,\omega)$, see Remark \ref{rem:truncation-Arthur-compatibility} for the notation. 
Suppose that $f(x^{-1}\gamma x)\neq 0$ for some $\gamma\in \tfo$ and 
\[
\sum_{\xi\in Q(\Q)\bs G(\Q)} F^Q(\xi x,T)\tau_Q(H_0(\xi x)-T)=1.
\]
Then, by using \cite{LW13}*{Lemma 3.5.5}, Lemma \ref{lem:fsa-compact-support}, $Q^+=G$, and $F^G(x,T)=0$, we can prove that
\begin{equation}\label{eq:R=Gbound}
 \|H^-\|\ll_\Omega 1   \quad \text{and} \quad   d(T)-c_\Omega \ll_\Omega \|H^\theta \| \ll_\Omega 1+\|T\|.
\end{equation}
for some constant $c_\Omega >0$. 
These support bounds on the domain of integration in \eqref{eq:R=G K_G} alone do not determine the asymptotic size of $M_{\tfo,Q,G}^T(f)$. 
The following example shows that such a contribution need not decay as $d(T)$ tends to infinity, even for a fixed compactly supported test function.
\begin{example}\label{ex:cga-gl3-nondecay}
Let $G=\GL_3$ and let $P_0$ be the minimal parabolic subgroup consisting of the upper triangular matrices.
Put
\[
 J\coloneqq \begin{pmatrix}0&0&1\\0&-1&0\\1&0&0\end{pmatrix},
 \qquad \theta(g)\coloneqq J\,{}^tg^{-1}J^{-1}.
\]
Let $G'=G\rtimes\langle\theta\rangle$ and $\tG=G\delta$,
where $\delta=1\rtimes \theta$.
Take $Q$ of type $(2,1)$ and $R=G$, then $Q^+=R^-=G$. 
Choose $K=\mathrm O(3)K_f$, where $K_f=\GL_3(\widehat\Z)$.
Let $\fo_\delta$ be the coarse class of $\delta$ and set
$T_t\coloneqq (t,0,-t)$, so that $\theta(T_t)=T_t$ and $d(T_t)=t$.
By using \eqref{eq:R=G}, we show that some fixed nonnegative $f\in C_c^\infty(\tG(\A)^1)$ satisfies
\begin{equation}\label{eq:gl3-nondecay-conclusion}
 M_{\fo_\delta,Q,G}^{T_t}(f)\gg_f1
 \qquad(\forall t\gg1).
\end{equation}

Write
\[
 n(p,q,r)\coloneqq\begin{pmatrix}1&p&r\\0&1&q\\0&0&1\end{pmatrix}.
\]
We have $N_{P_0}(\Q)\delta\subset\fo_\delta$.
Indeed, putting $z=(p+q)/2$, $b=(p-q)/2$, and
$u\coloneqq n(b,0,(r-z^2/2+b^2)/2)$, one checks that
\[
 n(p,q,r)\delta
 =u\bigl(n(z,z,z^2/2)\delta\bigr)u^{-1}.
\]
The unipotent element $n(z,z,z^2/2)$ is fixed by $\theta$,
so the semisimple part of the expression in parentheses is $\delta$.

Put
\[
 a_Q(X)\coloneqq\diag(e^X,e^X,e^{-2X}),\quad
 a_M(Y)\coloneqq\diag(e^Y,e^{-Y},1),\quad a_H\coloneqq a_M(Y)a_Q(X).
\]
In the standard $\GL_2$ cusp of $M_Q(\Q)\backslash M_Q(\A)^1$,
restrict the finite component to
$\GL_2(\widehat\Z)\times\widehat\Z^\times$ and the real
unipotent coordinate to a fixed interval.
For $Y\gg1$, the rank-one truncation condition is
$F^Q(m,T_t)=1$ precisely when $Y\le t/2$.
Also, $\tau_Q(H_0(a_Q(X))-T_t)=1$ precisely when $X>t/2$.
Thus both factors equal one on the rectangle
\begin{equation}\label{eq:gl3-nondecay-rectangle}
 X=t/2+u,\qquad Y=t/2-v,\qquad 1\le u,v\le2.
\end{equation}
This cusp parametrization has uniformly bounded multiplicity.
The measure in the integral, apart from compact variables of
fixed positive volume, is
\[
 e^{-6X}\,dX\;e^{-2Y}\,dY.
\]
Here $e^{-6X}$ comes from $\delta_Q(a_Q(X))^{-1}$ and
$e^{-2Y}$ from the Levi cusp measure.

Choose compact representatives for $N_Q(\Q)\backslash N_Q(\A)$,
with integral finite components.
Combining their real components with the Levi unipotent coordinate,
we may write the integration variable as $n_0a_Hk$, with $n_0$
in a fixed compact subset of $N_{P_0}(\A)$ and $k\in K$.
On \eqref{eq:gl3-nondecay-rectangle},
\[
 H=(X+Y,X-Y,-2X),\qquad
 a_H^{-1}\theta(a_H)=\diag(e^{u+v},e^{-2(u+v)},e^{u+v}).
\]
For integers satisfying
\begin{equation}\label{eq:gl3-nondecay-integer-box}
 |p|\le e^{2Y},\qquad |q|\le e^{3X-Y},\qquad |r|\le e^{3X+Y},
\end{equation}
all entries above the diagonal of $a_H^{-1}n(p,q,r)a_H$
have absolute value at most one.
Moreover, conjugation by $a_H^{-1}$ contracts the fixed compact
sets containing $n_0$ and $\theta(n_0)$.
Consequently, the pullbacks of
\[
 (n_0a_H)^{-1}n(p,q,r)\delta(n_0a_H)
\]
range over a fixed compact subset of $G(\A)^1$, with finite
components in $K_f$.
Choose a nonnegative compactly supported $f$ whose pullback is
bi-$K$-invariant and at least one on this compact set.
Such an $f$ can be chosen independently of $t$; since
$\theta(K)=K$, it satisfies $f_K=f$.

All the integer triples in \eqref{eq:gl3-nondecay-integer-box}
therefore contribute at least one to the sum over $\fo_\delta$.
Their number is
\[
 \gg e^{2Y}e^{3X-Y}e^{3X+Y}=e^{6X+2Y}.
\]
Restricting the nonnegative integral to
\eqref{eq:gl3-nondecay-rectangle}, we obtain
\[
 M_{\fo_\delta,Q,G}^{T_t}(f)
 \gg \int_1^2\int_1^2 e^{-6X-2Y}e^{6X+2Y}\,du\,dv=1.
\]
This proves \eqref{eq:gl3-nondecay-conclusion}.
\end{example}

\begin{remark}\label{rem:TW}
In \cite{TW26}*{Theorem 3.2}, we proved an asymptotic formula for twisted traces by using the twisted trace formula. 
As explained above, \cite{Par19}*{Theorem 4.1 (3)} contains an error, and the correct formulation for compactly supported test functions is given in Theorem \ref{thm:cga}.
Since \cite{TW26}*{Lemma 4.25} relies on \cite{Par19}*{Theorem 4.1 (3)}, minor modifications to the proof of \cite{TW26}*{Theorem 3.2} are needed.
We replace the use of \cite{Par19}*{Theorem 4.1 (3)} by Theorem \ref{thm:cga}, leaving the rest of the argument unchanged.
We briefly explain how to make these modifications here.    

To correct the argument, it suffices, by the cuspidality of $h_{v_0}$ (see \cite{TW26}*{(4.12)}), to replace $F^G(g,T)$ in the integrals defining $I^T_{\mathrm{negl},1}(h_j)$ and $I^T_{\mathrm{negl},2}(h_j)$ in \cite{TW26}*{p.44} by
\begin{equation}\label{eq:TWcorrection}
F^G(g,T)+\sum_{\substack{Q\subsetneq G \\ Q^+=G}}\sum_{\xi\in Q(\Q)\bs G(\Q)}F^Q(\xi g,T)\, \tau_Q(H_0(\xi g)-T).
\end{equation}
With this replacement, Theorem \ref{thm:cga} gives the required bound for the contributions of the full coarse classes. The contribution of $\fo_{\mathrm{main}}\sqcup\fo_{\mathrm{negl,ss}}$ is estimated as in the proof of \cite{TW26}*{Lemma 4.25}, since the complementary truncation involves only parabolic subgroups $Q$ with $Q^+\subsetneq G$.
Furthermore, $F^G(g,T)$ must also be replaced by \eqref{eq:TWcorrection} in \cite{TW26}*{Lemma 4.27}, but \eqref{eq:R=Gbound} ensures that its proof remains valid without further changes.
These modifications complete the correction to the proof of \cite{TW26}*{Theorem 3.2}.
\end{remark}

\section{On the continuity of the coarse geometric expansion}
\label{gr:section}

We retain the notation and measures of \S~\ref{sec:preliminaries} and \S~\ref{sec:cga}.
Although some of the notation in \S~\ref{sec:preliminaries} and \S~\ref{sec:cga}, such as $K_{\tP}(x,y)$ and $B^T_{\fo,Q,R}(f)$, was originally introduced for $f\in C^\inf_c(\tG(\A)^1)$ or $f\in C^\inf_\Omega(\tG(\A)^1;K_0)$, we use the same definitions for $f\in \cC(\tG(\A)^1;K_0)$ as well.

\subsection{Absolute convergence theorem}

The following is the main result of this section.
\begin{theorem}\label{thm:main theorem 1}
There exist constants $C$, $d_*>0$ and integers $m$, $s$, $r\geq0$,
with $d_*>d_0$, such that the following assertions hold for
every open compact subgroup $K_0\subset G(\A_f)$ and 
every $f\in\cC(\tG(\A)^1;K_0)$.  
\begin{enumerate}[(1)]
\item For any $T\in\fa_{M_0}^G$ satisfying $d(T)>d_*$,
\[
 \sum_{\tfo\in\tcO}\int_{G(\Q)\bs G(\A)^1}|k_{\tfo}^T(x,f)| \,dx
\leq C\level(K_0)^s
       (1+\|T\|)^r\|f\|_{\tG(\A)^1,m}.
\]
Hence, the integrals $J^{T,\tG}(f)$ and $J_{\tfo}^{T,\tG}(f)$ are absolutely convergent, and we have
\[
J^{T,\tG}(f)=\sum_{\tfo\in\tcO}J_{\tfo}^{T,\tG}(f).
\]

\item 
$J^{T,\tG}(f)$ and $J_{\tfo}^{T,\tG}(f)$ are polynomials in
$T$ of degree at most
$\dim\fa_{M_0}^{G}$.
Their coefficients are continuous linear forms in $f$.
For a fixed monomial basis, the sum over $\tfo$ of the absolute
values of each coefficient is bounded by $C\level(K_0)^s\|f\|_{\tG(\A)^1,m}$, after increasing $C$ if necessary.

\item For any $T\in\fa_{M_0}^G$ satisfying $d(T)>d_*$, the approximation estimate
\begin{align*}
 \sum_{\tfo\in\tcO}
   |J_{\tfo}^{T,\tG}(f)-M_{\tfo}^T(f)|
 &\leq \sum_{\tfo\in\tcO}\int_{G(\Q)\bs G(\A)^1} |k_{\tfo}^T(x,f)-m_{\tfo}^T(x,f)|\, dx  \nonumber\\
 &\leq C\level(K_0)^s
       (1+\|T\|)^r e^{-d(T)}\|f\|_{\tG(\A)^1,m}
\end{align*}
holds. 
Note that $\sum_{\tfo\in\tcO}\int_{G(\Q)\bs G(\A)^1}|m_{\tfo}^T(x,f)| \,dx\le C\level(K_0)^s\, (1+\|T\|)^r\, \|f\|_{\tG(\A)^1,m}$. 
\end{enumerate}
\end{theorem}
In the remainder of this section, we will explain the proof of this theorem.

\subsection{Reduction of the proof of Theorem \ref{thm:main theorem 1} (1) (3) to Theorem \ref{thm:bound 1}}\label{sec:reduction(1)(3)}
We first show that the kernel functions appearing in the statement of the theorem are well defined at every point and continuous in $\cC(\tG(\A)^1;K_0)$.
\begin{lemma}\label{lem:sup kernel}
Fix a compact subset $D$ of $G(\A)^1$. 
For any $f\in\cC(\tG(\A)^1;K_0)$, there
are $c_D>0$ and integers $m$, $s_D\geq0$ such that
\[
 \sup_{x\in D}\sum_{\gamma\in \tM_P(\Q)}
       \int_{N_P(\A)}
       |f(x^{-1}\gamma nx)|\,dn
 \le c_D\, \level(K_0)^{s_D} \,
                     \|f\|_{\tG(\A)^1,m}.
\]
For fixed $T$ with $d(T)>d_0$, the kernels $k_{\tfo}^T$, $B_{\tfo,Q,R}^T$, and $m_{\tfo}^T$
are therefore well defined by their original formulas.
\end{lemma}
\begin{proof}
Let $\cN'$ denote a compact fundamental subset of $N_P(\Q)\bs N_P(\A)$. 
By the Fubini--Tonelli theorem,
\begin{align*}
\sum_{\gamma\in \tM_P(\Q)} \int_{N_P(\A)} |f(x^{-1}\gamma nx)|\,dn
       = &\int_{N_P(\Q)\bs N_P(\A)} \sum_{\gamma\in \tP(\Q)}|f(x^{-1}\gamma nx)|\,dn \\
       \le & \int_{\cN'} \sum_{\gamma\in \tG(\Q)}|f(x^{-1}\gamma nx)|\,dn.    
\end{align*}
Set $r\coloneqq \level(K_0)$. 
By Lemma \ref{lem:levelbound}, there exists a natural number $q(\cN'D)$ such that the function $\tG(\A)^1\ni g\mapsto f(x^{-1}g nx)$ is in $\cC(\tG(\A)^1;K(r\, q(\cN'D)))$. 
By Lemma \ref{lem:1}, there exist integers $m'$, $s'\ge 0$ such that
\[
\sum_{\gamma\in \tG(\Q)}|f(x^{-1}\gamma nx)|\ll_D \level(K(r\, q(\cN'D)))^{s'} \|f\|_{\tG(\A)^1,m'}
\]
holds for any $x\in D$ and any $n\in \cN'$. 
Since $\level(K(r\,q(\cN'D)))\le q(\cN'D)\, \level(K_0)$, these two inequalities imply the first assertion of the lemma.

By \cite{LW13}*{Lemma~3.7.2}, for fixed $T$, the factor
$\widehat\tau_{\tP}(H_P(\xi x)-T)$ selects only finitely
many classes in $P(\Q)\backslash G(\Q)$ uniformly
in $x\in D$.
Hence, the second assertion for $ k_{\tfo}^T$ holds. 
Since $F^Q(\xi x,T)\tau_Q^P(H_P(\xi x)-T)$ also only finitely
many classes in $Q(\Q)\backslash G(\Q)$ uniformly
in $x\in D$, we obtain the same assertion for $B_{\tfo,Q,R}^T$ and $m_{\tfo}^T$ by \eqref{eq:lw911-sigma-identity}.
\end{proof}

Take a test function $f\in \cC(\tG(\A)^1;K_0)$ and let $Q$, $R\in\cF^G(M_0)$ such that $P_0\subset Q\subset R$ and $Q^+\subset R^-$. 
By Lemma \ref{lem:sup kernel}, as in the proof of Theorem \ref{thm:cga}, 
it suffices to obtain an upper bound for the following integral $\cR_{Q,R}(T,f)$. 
\[
 \cR_{Q,R}(T,f)
  \coloneqq\sum_{\tfo\in\tcO}\int_{Q(\Q)\bs G(\A)^1}
      F^Q(x,T)\widetilde\sigma_Q^R(H_0(x)-T)
      \, | k_{\tfo,Q,R}(x,f)|\,dx,
\]
\[
 k_{\tfo,Q,R}(x,f)
   \coloneqq \sum_{\tQ^+\subset \tP\subset \tR^-}
        (-1)^{\dim\fa_{\tP}^{\tG}} K_{\tP,\tfo}(x,x).
\]

For $\tS\in\cF^{\tG}(\tM_0)$ with $Q^+\subset S$, set
\[
 \tM_S(\Q)_{Q}^{\circ}
   \coloneqq \tM_S(\Q)\setminus
       \bigcup_{\tQ^+\subset \tP\subsetneq \tS}\tP(\Q).
\]
Note that $\gamma^{-1}\tM_S(\Q)_{Q}^{\circ}\gamma=\tM_S(\Q)_{Q}^{\circ}$ holds for any $\gamma\in M_{S}(\Q)\cap Q(\Q)$. 
Since 
\[
k_{\tfo,Q,R}(x,f)=\sum_{Q^+\subset S\subset R^-}\cE_{Q,R;S}(x,\tM_S(\Q)_{Q}^{\circ}\cap\tfo,f),
\]
one can apply Proposition \ref{prop:tfl47} to $k_{\tfo,Q,R}(x,f)$. 
We note that $\cE_{Q,R;S}(x,\tM_S(\Q)_{Q}^{\circ}\cap\tfo,f)$ is denoted by 
$k_{\tfo,1,2,3}$ on \cite{FL16}*{p.444} and $k_{\mathscr{O},Q,R,S}(x)$ on \cite{Par19}*{p.~545}.
Therefore, 
there exist $s\in\Z_{>0}$ and $X_1,\dots,X_m\in\cU(\fg_\C^1)$, which are independent of $R$ and $Q$, such that
\begin{align*}
\cR_{Q,R}(T,f)\ll& \level(K_0)^s\sum_{\tfo\in\tcO} \sum_{\tQ^+\subset\tS \subset \tR^-}\sum_{\gamma\in \tM_S(\Q)_{Q}^{\circ}\cap\tfo} \int_{Q(\Q)\bs Q(\A)\cap G(\A)^1}dq \int_K \, dk  \\  
      & F^Q(q,T)\widetilde\sigma_Q^R(H_0(q)-T)  e^{-\Gamma_{Q,R;S}(T,H_0(q))} \int_{N_{S}(\A)} \sum_{j=1}^m|f_j(k^{-1}q^{-1}\gamma nqk)|\,dn,
\end{align*}
where $f_j\coloneqq f*X_j$.
By Lemma \ref{lem:1} (3), there exists a right-$K$-invariant function $\tf\in \cC(\tG(\A)^1)$ such that 
\[
\sum_{j=1}^m|f_j(k^{-1}q^{-1}\gamma nqk)|\le \tf(k^{-1}q^{-1}\gamma nq)
\]
and \eqref{eq:3} holds for $f$ and $\widetilde f$, with the derivative order on the right increased by $\max_j\deg X_j$.
Hence, 
\begin{align*}
\cR_{Q,R}(T,f)\ll& \level(K_0)^s \sum_{\tQ^+\subset\tS \subset \tR^-}\sum_{\gamma\in \tM_S(\Q)_{Q}^{\circ}} \int_{Q^{M_S}(\Q)\bs M_S(\A)\cap G(\A)^1}dx   \\  
      & F^Q(x,T)\widetilde\sigma_Q^R(H_0(x)-T)  e^{-\Gamma_{Q,R;S}(T,H_0(x))}\, \int_{N_{S}(\A)}\tf_{K}(x^{-1}\gamma xn)\, dn \nonumber
\end{align*}
where $Q^{M_S}\coloneqq Q\cap M_S$ and $\tf_{K}(x)\coloneqq \int_{K}\tf(k^{-1}x)\, dk$. 
Note that $\tf_{K}$ is bi-$K$-invariant. 

For $P_1$, $P_2\in\cF^G(M_0)$ with $P_0\subset P_1\subset P_2$, 
we set
\[
\cJ_{P_1,P_2}(T)\coloneqq\begin{cases}
    e^{-d(T)} & \text{if $\tP_1^+\subsetneq\tP_2^-$,}\\
    1 & \text{if $\tP_1^+=\tP_2^-$.}
\end{cases}
\]
For $P_1\subset P_2\subset P_3$ and $H_{P_1}^{P_3}\in \fa_{P_1}^{P_3}$, we write $H_{P_1}^{P_3}=H_{P_1}^{P_2}+H_{P_2}^{P_3}$ for the decomposition of $H_{P_1}^{P_3}$ with respect to $\fa_{P_1}^{P_3}=\fa_{P_1}^{P_2}\oplus\fa_{P_2}^{P_3}$.
The following lemma is a twisted analogue of \cite{FL16}*{Lemma 4.8}
\begin{lemma}\label{lem:analogueFL48}
Let $y\in M_S(\A)^1$, $T\in\fa_{M_0}^G$, $H_S^\theta\in\fa_{\tS}^{\tG}$, and $H_S^-\in(\fa_S^G)^{1-\theta}$.
If $d(T)>d_0$, then
\begin{align*}
&\int_{\fa_{\tS}^{\tG}}
 \widetilde\sigma_Q^R\bigl(H_0(y)+H^-_S+H^\theta_S-T\bigr)
 e^{-\Gamma_{Q,R;S}(T,H_0(y)+H^-_S+H^\theta_S)}\,dH^\theta_S\\
&\ll
\cJ_{S,R}(T)
\tau_Q^S(H_0(y)-T)
\left(1+\|H_Q(y)^S-T_Q^S\|+
 \|H^-_S\| +\|T\| \right)^{\dim\fa_{\tR^-}^{\tG}}.
\end{align*}
\end{lemma}
\begin{proof}
For simplicity, we put
\[
Y\coloneqq H_Q(y)^S-T_Q^S\in \fa_Q^S,
\quad
Z\coloneqq H_S^--(T_S^{G})^-\in(\fa_S^G)^{1-\theta},
\]
\[
d_R^G\coloneqq \dim\fa_{\tR^-}^{\tG},
\quad
d_S^R\coloneqq \dim\fa_{\tS}^{\tR^-},
\]
\[
I(T,y,H_S^-)
\coloneqq 
\int_{\fa_{\tS}^{\tG}}
\tsigma_Q^R(H_0(y)+H_S^-+H_S^\theta-T)\,
e^{-\Gamma_{Q,R;S}(T,H_0(y)+H_S^-+H_S^\theta)}\,dH_S^\theta.
\]
Since $\tsigma_Q^R(H_0(y)+H_S^-+H_S^\theta-T)=\tsigma_Q^R(H_Q(y)+H_S^-+H_S^\theta-T_Q)$, we have
\[
I(T,y,H_S^-)=\int_{\fa_{\tS}^{\tG}}
\tsigma_Q^R(Y+Z+H_S^\theta)\,
e^{-\Gamma_{Q,R;S}(T,Y+Z+H_S^\theta+T_Q)}\,dH_S^\theta,
\]
by change of variable $H_S^\theta\mapsto H_S^\theta+(T_S^G)^\theta$.
By $H_S^\theta=(H_S^{R^-})^\theta+H_{R^-}^\theta$,
Lemma \ref{lem:cga-cone-fibers} gives
\[
\widetilde\sigma_Q^R(Y+Z+H_S^\theta)=1
\quad\Longrightarrow\quad
\|H_{R^-}^\theta \|\ll\|Y+Z+(H_S^{R^-})^\theta\|.
\]
For each $\tbeta\in\Delta_{\tP_0}^{\tR^-}\setminus\Delta_{\tP_0}^{\tS}$, we set
\[
u_{\tbeta}\coloneqq \widetilde\beta((H_S^{R^-})^\theta)+\min_{\substack{\alpha\in\Delta_{P_0}\\
                 \widetilde\alpha=\widetilde\beta}} \alpha_Q(Y+Z).
\]
Since $\alpha_Q(H_{R^-}^\theta)=0$ for $\alpha\in \Delta_0^{R^-}$, we have 
\[
\Gamma_{Q,R;S}(T,Y+Z+H_S^\theta+T_Q)= \Gamma_{Q,R;S}(T,Y+Z+(H_S^{R^-})^\theta+T_Q)
 \ge  d_S^R \, d(T)+\sum_{\tbeta\in\Delta_{\tP_0}^{\tR^-}\setminus\Delta_{\tP_0}^{\tS}}u_{\tbeta} .  
\]
Since $\{ \talpha_{S} \mid \talpha\in \Delta_{\tP_0}^{\tR^-}\setminus\Delta_{\tP_0}^{\tS}\}$ forms a basis of $(\fa_{\tS}^{R^-})^*$, by $|\alpha_Q(Y+Z)|\ll \|Y\|+\|Z\|$, we obtain
\[
\|(H_S^{R^-})^\theta\| \ll \|Y\|+\|Z\|
+\sum_{\tbeta\in\Delta_{\tP_0}^{\tR^-}\setminus\Delta_{\tP_0}^{\tS}}u_{\tbeta}.
\]
Therefore, by $d(T)>d_0$,
\begin{align*}
I(T,y,H_S^-)
&\ll
\tau_Q^S(H_0(y)-T) \cJ_{S,R}(T)
\int_{(\fa_{\tS}^{\tR^-})_{y,H_S^-,+}}
e^{-\sum_{\widetilde\beta}u_{\widetilde\beta}}
\Bigl(
1+\|Y\|+\|Z\|
+\sum_{\widetilde\beta}u_{\widetilde\beta}
\Bigr)^{d_R^G} \, d(H_S^{R^-})^\theta \\ 
&\ll
\tau_Q^S(H_0(y)-T)\,\cJ_{S,R}(T)
\bigl(1+\|Y\|+\|Z\|\bigr)^{d_R^G},    
\end{align*}
where $(\fa_{\tS}^{\tR^-})_{y,H_S^-,+}\coloneqq \{ (H_S^{R^-})^\theta\in \fa_{\tS}^{\tR^-} \mid  u_{\tbeta}>0 \quad (\forall \tbeta\in \Delta_{\tP_0}^{\tR^-}\setminus\Delta_{\tP_0}^{\tS} ) \}$.
This completes the proof.
\end{proof}
By this lemma, putting $\ell \coloneqq \dim\fa_{\tR^-}^{\tG}$, 
\begin{align*}
\cR_{Q,R}(T,f)\ll& \level(K_0)^s \,  (1+\|T\|)^{\ell} \sum_{\tQ^+\subset\tS \subset \tR^-}  \cJ_{S,R}(T)\sum_{\gamma\in \tM_S(\Q)_{Q}^{\circ}} \int_{Q^{M_S}(\Q)\bs M_S(\A)^1}dx   \\  
      & F^Q(x,T)\, \tau_Q^S(H_0(x)-T) \, (1+\|H_Q(x)^S-T_Q^S\|)^{\ell}\, \tf_{S,\ell}(x^{-1}\gamma x). \nonumber
\end{align*}
Here, $\tf_{S,\ell}(x)=\int_{N_S(\A)}\int_{(A_S^G)^{1-\theta}}(1+ \| H_S(a)\| )^{\ell}\tf_{K}(axn)\, da \, dn$. 
By Proposition \ref{prop:biKP} and Lemma \ref{lem:1} (3),
\[
\tf_{S,\ell}\in \cC(\tM_S(\A)^1), \quad  \|\tf_{S,\ell}\|_{\tM_S(\A)^1,m} \ll \|\tf\|_{\tG(\A)^1,m+\dim(\fa_S^G)^\theta} \ll \level(K_0)^{s'}\,  \| f \|_{\tG(\A)^1,m+m'}
\]
for some $s'>0$ and $m'>0$. 
In addition, $\tf_{S,\ell}$ is bi-$K_{M_S}$-invariant and non-negative, where $K_{M_S}\coloneqq K\cap M_S(\A)$. 
Consequently, the assertions (1) and (3) of Theorem \ref{thm:main theorem 1} follow from Theorem \ref{thm:bound 1}, since $K_{M_S}$ is also admissible relative to $M_0$ in $M_S$, see \cite{Art81}*{p.9}. 
\begin{remark}
    Note that Parab's paper \cite{Par19} contains no discussion whatsoever of either Lemma \ref{lem:analogueFL48} or Proposition \ref{prop:biKP}.
\end{remark}

\subsection{Statement and Proof of Theorem \ref{thm:bound 1}}

\begin{theorem}\label{thm:bound 1} 
    Fix a natural number $l\ge 0$. 
    Then there exist integers $r\geq0$ and $s>0$ such that
\begin{align}\label{eq:thm mellin}
 &\int_{P(\mathbb Q)\backslash G(\mathbb A)^1}
 F^P(g,T)\tau_P(H_P(g)-T)
 (1+\|H_P(g)-T_P\|)^l
 \sum_{\gamma\in \tG(\Q)^\circ_{P}}
 f(g^{-1}\gamma g) \, dg \\
 &\qquad\ll_l
 (1+\|T\|)^r \cJ_{P,G}(T) \, \| f\|_{\tG(\A)^1,s}  \nonumber
\end{align}
for any $P\in\cF^G(M_0)$ with $P_0\subset P$, any bi-$K$-invariant and non-negative function $f\in \cC(\tG(\A)^1)$, and any $T\in\mathfrak a_{P_0}^G$ such that $\langle\alpha,T-T_1\rangle\geq0$ for all $\alpha\in\Delta_{P_0}$.
\end{theorem}
\begin{proof}
We first suppose that $f$ is compactly supported. 
Put $f_\delta(x)\coloneqq f(x\delta)$.
Since $f$ is bi-$K$-invariant, $f_\delta$ is left-$K$-invariant and right-$\theta(K)$-invariant. 
Set
\[
\cD_P(T)\coloneqq \{H\in\fa\mid 
\tau_{P_0}^{P}(H-T_1)=1,\quad
\widehat\tau_{P_0}^{P}(T-H)=1, \quad \tau_P(H-T)=1\}.
\]
Then, the left hand side of \eqref{eq:thm mellin} is bounded by a constant multiple of 
\begin{equation}\label{eq:r1}
\int_{\cD_P(T)}\int_{P_0(\Q)\backslash P_0(\A)^1}
 (1+\|H_P-T_P\|)^l \, e^{-2\rho_{P_0}(H)}
 \sum_{\gamma\delta\in \tG(\Q)^\circ_{P}}
 f_\delta(a^{-1}p^{-1}\gamma \theta(pa)) \, dp\, dH, \quad a=\exp(H).
\end{equation}

For simplicity, we put
\[
\fa\coloneqq \fa_{M_0}^{G},\quad W\coloneqq W_0^G
,\quad \Delta\coloneqq \Delta_{P_0}, \quad N_0\coloneqq N_{P_0}.
\]
Choose a representative $n_w\in N_G(M_0)(\Q)$ for each $w\in W$. 
For $w\in W$, let $I(w)\subset\Delta$ be the support of a reduced expression for $w$, and put
\[
N_w\coloneqq N_0\cap n_wN_0n_w^{-1}.
\]
Let $Q(w)$ be the smallest standard parabolic subgroup containing $n_w$. Then, $I(w)=\Delta_{P_0}^{Q(w)}$.
For $U\subset \Delta$, we write $\sat_\theta(U)$ for the smallest $\theta$-stable subset containing $U$. 
Then, 
\begin{equation}\label{eq:r2}
\tG(\Q)_P^\circ=\bigsqcup_{\substack{w\in W \\ \sat_\theta(\Delta_{P_0}^P\cup I(w))=\Delta }} P_0(\Q) n_w P_0(\Q) \delta.
\end{equation}
For $\mathbf{f}\in C_c^\inf(G(\A)^1)$, we define
\[
F_{\mathbf{f}}(g)\coloneqq \int_{P_0(\A)^1}\mathbf{f}(p\, g)\, dp, \quad g\in G(\A)^1,
\]
\begin{equation*}
J_{P,w,l}^{T}(\mathbf{f})\coloneqq \int_{\cD_P(T)}\int_{N_{w^{-1}}(\A)\backslash N_0(\A)}
(1+\|H_P-T_P\|)^l  F_{\mathbf{f}}( \exp(-L_w(H)) n_wn) \, e^{2\rho_{w^{-1}}(\theta H)}\, dn \,dH, 
\end{equation*}
where $L_w\coloneqq 1-w\theta$ and $[e^H\mapsto e^{2\rho_w(H)}]$ is the modulus function of $A_0^G$ on $N_w(\A)\bs N_0(\A)$.
By \cite{FL16}*{Proof of Theorem 3.1}, \cite{FL16}*{Lemma 3.3}, \eqref{eq:r1} and \eqref{eq:r2}, there are $s'>0$ and $X_1,\dots,X_{m'}\in \cU(\fg^1_\C)_{\le s'}$ such that the left hand side of \eqref{eq:thm mellin} is bounded by
\begin{equation*}
\sum_{\substack{w\in W \\ \sat_\theta(\Delta_{P_0}^P\cup I(w))=\Delta }} \sum_{j=1}^{m'} J_{P,w,l}^{T}(|f_{\delta,j}|)    
\end{equation*}
where $f_{\delta,j}\coloneqq f_\delta*X_j$ and $|f_{\delta,j}|(x)\coloneqq |f_{\delta,j}(x)|$, see also \cite{FM21}*{Proof of Lemma 4.6}. 
In Lemma \ref{lem:1} (3), we may replace $K$ by $\theta^{-1}(K)$.
Hence, by Lemma \ref{lem:1} (3), there exists a bi-$K$-invariant function $\tf\in C_c^\inf(G(\A)^1)$ such that, for any $m\in\Z_{\ge 0}$,
\begin{align*}
& \tf(x)\ge|f_{\delta,j}(x)| \;\; (\forall x\in G(\A)^1,\;\; 1\le j\le m'), \\
&  \|\tf(x)\|_{G(\A)^1,m} \ll_m  \|f(x)\|_{\tG(\A)^1,m+s'+\dim G-\mathrm{rank}X^*(G)}   .  
\end{align*}
Note that this $\tf$ differs from the $\tf$ used in \S\ref{sec:reduction(1)(3)}.
Thus, the left hand side of \eqref{eq:thm mellin} is bounded by a constant multiple of
\begin{equation}\label{eq:r3}
\sum_{\substack{w\in W \\ \sat_\theta(\Delta_{P_0}^P\cup I(w))=\Delta }} J_{P,w,l}^{T}(\tf)    
\end{equation}
After proving several lemmas, we establish the estimate $J_{P,w,l}^{T}(\tf) \ll (1+\|T\|)^r \cJ_{P,G} \|\tf\|_{G(\A)^1,s''}$ in Lemma \ref{lem:th45-bruhat} for some $s''>0$. 
This completes the proof of the theorem for $f\in C^\inf_c(\tG(\A)^1)$.
After that, we use Fatou's lemma to extend the theorem to $\cC(\tG(\A)^1)$.

Set $\rho\coloneqq \rho_{P_0}$.
For $\mathbf{f}\in C_c^\inf(G(\A)^1)$, define
\[
 \varphi_{\mathbf{f}}(g,\lambda)
 \coloneqq \int_\fa F_{\mathbf{f}}(\exp(H) g) \, e^{-(\lambda+\rho)(H)}\,dH, \qquad \lambda\in\fa_\C^*.
\]
We also set $ \varphi_{\mathbf{f}}(\lambda)\coloneqq \varphi_{\mathbf{f}}(1,\lambda)$. 
If $\mathbf{f}$ is bi-$K$-invariant, then $\varphi_{\mathbf{f}}$ is $W$-invariant. 
Let $\sC_\rho$ denote the convex hull of the points $W\rho$ in $\fa^*$.
For $\sigma\in\mathcal \sC_\rho$, put
\[
 \Psi_{\mathbf{f},\sigma}(H)\coloneqq e^{(-\rho-\sigma)(H)}F_{\mathbf{f}}(\exp(H)).
\]
Set $n\coloneqq \dim\fa$, and fix a coordinate on $\fa$. 
For $\mathbf{a}=(a_1,\dots,a_n)\in(\Z_{\ge0})^n$, $H=(H_1,\dots,H_n)$,
\[
|\mathbf{a}|=a_1+\cdots+a_n, \quad  \partial^{\mathbf{a}}_H=\left(\frac{\partial}{\partial H_1}\right)^{a_1}\cdots\left(\frac{\partial}{\partial H_n}\right)^{a_n}.
\]
\begin{lemma}
\label{lem:th45-abel}
Suppose that $\mathbf{f}$ is bi-$K$-invariant, hence $\varphi_{\mathbf{f}}$ is $W$-invariant. 
Then, for every integer \(b\geq0\),
\begin{equation*}
 \sum_{|\mathbf a|\leq b}
 \|\partial_H^{\fa}\Psi_{\mathbf{f},\sigma}\|_{L^1(\fa)}
 \ll_b \|\mathbf{f}\|_{G(\A)^1,b},
 \qquad \sigma\in\sC_\rho.
\end{equation*}
In particular, this estimate is uniform up to the boundary of
$\sC_\rho$.
\end{lemma}
\begin{proof}
First suppose $\mathbf{f}\ge 0$. 
If $\sum_{w\in W}t_w=1$ $(t_w>0)$ and $\sigma=\sum_{w\in W}t_w(w\rho)\in\sC_\rho$, then
\[
\varphi_{\mathbf{f}}(\sigma)=\int_{G(\A)^1} \prod_{w\in W} \left(\mathbf{f}(g)e^{(\rho-w\rho)(H)}\right)^{t_w}dg\le \prod_{w\in W}\varphi_{\mathbf{f}}(w\rho)^{t_w}=\|\mathbf{f}\|_{L^1(G(\A)^1)}
\]
by H\"older's inequality and $\varphi_{\mathbf{f}}(-\rho)=\|\mathbf{f}\|_{L^1(G(\A)^1)}$.
Hence, we have $\varphi_{\mathbf{f}}(\sigma)\le \|\mathbf{f}\|_{L^1(G(\A)^1)}$ for any $\sigma\in\sC_\rho$. 
For $X\in\fa$, we have $\partial_X\Psi_{\mathbf{f},\sigma} =\Psi_{\mathbf{f}*X,\sigma}+(-\rho-\sigma)(X)\Psi_{\mathbf{f},\sigma}$. 
Thus $\partial_H^{\fa}\Psi_{\mathbf{f},\sigma}$ is a linear combination of the functions $\Psi_{\mathbf{f}*H^{\fb},\sigma}$ $(|\fb|\le |\fa|)$, with bounded coefficients and at most $2^{|\fa|}$ terms. 
As in the proof of Lemma \ref{lem:biK-absolute-global}, the absolute values of the finitely many derivatives occurring above admit a common nonnegative bi-\(K\)-invariant majorant $Q_b$ satisfying $\|Q_b\|_{L^1(G(\A)^1)}\ll_b\|\mathbf{f}\|_{G(\mathbb A)^1,b}$. 
Applying the preceding Hölder estimate to $Q_b$ proves the assertion.
\end{proof}

Since our aim is to bound \eqref{eq:r3} from above, we assume from now on that $w\in W$ satisfies 
\[
\operatorname{sat}_\theta\bigl(\Delta_{P_0}^P\cup I(w)\bigr) =\Delta.
\]
Set
\[
\fa^*_{++,w}\coloneqq \{ \lambda=\sum_{\alpha\in\Delta}a_\alpha\alpha\in\fa^* \mid a_\alpha\geq0 \quad (\forall \alpha\in\Delta),\quad \langle\lambda,\beta^\vee\rangle=0 \quad (\forall \beta\in I(w))\},
\]
\[
 \eta_w(\lambda)\coloneqq \lambda-\theta^{-1}w^{-1}\lambda.
\]
Recall that $\widehat\Delta^\vee\coloneqq \widehat\Delta^\vee_{P_0}=\{ \varpi_\alpha^\vee \mid \alpha\in\Delta \}$, which is the basis of $\fa$ dual to $\Delta$.
\begin{lemma}
\label{lem:th45-parameter}
There exist $\chi\in\fa^*_{++,w}$ and $\varepsilon_0>0$, $c>0$ such that, for
\begin{equation}\label{eq:th45-parameter}
 \sigma_T\coloneqq \rho-\frac{\varepsilon_0}{1+\|T\|}\chi,
\end{equation}
we have
\begin{equation}\label{eq:th45-parameter-properties}
 \sigma_T\in\sC_\rho,
 \qquad
 \eta_w(\sigma_T)(\varpi_\alpha^\vee)\geq\frac{c}{1+\|T\|}
 \quad\bigl(\alpha\in\Delta\setminus\Delta_{P_0}^P\bigr).
\end{equation}
\end{lemma}
\begin{proof}
Put 
\[
 J(w)\coloneqq \{\alpha\in\Delta\setminus\Delta_{P_0}^P\mid 
                         \theta\alpha\notin I(w)\}.
\]
Recall that $Q(w)$ is the smallest standard parabolic subgroup containing $n_w$ and $I(w)=\Delta_{P_0}^{Q(w)}$.
Let $\rho_{I(w)}\coloneqq \rho_{P_0\cap M_{Q(w)}}$.
Then 
\[
\chi_0\coloneqq \rho-\rho_{I(w)}=\sum_{\alpha\in\Delta}a_\alpha^{(0)}\alpha
\]
is $W_0^{Q(w)}$-invariant and $a_\alpha^0\ge0$.
For $\alpha\in J(w)$, the difference
$a_{\theta\alpha}^{(0)}-a_\alpha^{(0)}$ is the coefficient of $\alpha$ in $\rho_{I(w)}$, since
\[
\rho_{I(w)}-\theta^{-1}\rho_{I(w)}=\theta^{-1}\chi_0-\chi_0=\sum_{\alpha\in\Delta}(a_{\theta\alpha}^{(0)}-a^{(0)}_\alpha)\alpha.
\]
It is positive on $J(w)\cap I(w)$ and zero on $J(w)\setminus I(w)$.

On the vertices $\Delta\setminus I(w)$, draw an edge $\alpha\longrightarrow\theta\alpha$ for each
$\alpha\in J(w)\setminus I(w)$.
If there exists a directed cycle ($\theta$-orbit) disjoint from
$\Delta_{P_0}^P\cup I(w)$, then it contradicts $\operatorname{sat}_\theta\bigl(\Delta_{P_0}^P\cup I(w)\bigr)
 =\Delta$.
Hence, every directed cycle ($\theta$-orbit) has an intersection with $\Delta_{P_0}^P\cup I(w)$. 
Choose nonnegative numbers $\{ b_\alpha \}_{\alpha\in\Delta\setminus I(w)}$ strictly increasing along these edges.
We can extend $\{ b_\alpha \}_{\alpha\in\Delta\setminus I(w)}$ uniquely to $\{ b_\alpha \}_{\alpha\in\Delta}$ by
\[\chi_1\coloneqq \sum_{\alpha\in\Delta}b_\alpha\alpha, \qquad 
 \langle  \chi_1,    \beta^\vee\rangle=0
 \qquad(\forall \beta\in I(w)).
\]
The resulting coefficients $\{ b_\alpha \}_{\alpha\in I(w)}$ are nonnegative, because the Cartan submatrix on
$I(w)$ has nonnegative inverse, while its off-diagonal entries are
nonpositive.
For a sufficiently small constant $u>0$, we set
\[
 \chi\coloneqq \chi_0+u\chi_1=\sum_{\alpha\in\Delta} a_\alpha \alpha,\quad a_\alpha=a_\alpha^{(0)}+ub_\alpha.
\]
By the above arguments, $\chi$ annihilates the coroots in $I(w)$ and satisfies
$a_{\theta\alpha}-a_\alpha>0$ for every $\alpha\in J(w)$.
In particular, $w\chi=\chi$ holds.

Since \(\theta\rho=\rho\), the usual inversion-root identity gives
\begin{equation}\label{eq:th45-eta-rho}
 \eta_w(\rho)=\theta^{-1}(\rho-w^{-1}\rho)
 =\sum_{\alpha\in\theta^{-1}I(w)}c_\alpha\alpha,
 \qquad c_\alpha\geq1.
\end{equation}
Then,
\[
 \eta_w(\sigma_T)(\varpi_\alpha^\vee)
 =\eta_w(\rho)(\varpi_\alpha^\vee)
   +\frac{\varepsilon_0}{1+\|T\|}
                      (a_{\theta\alpha}-a_\alpha).
\]
For $\alpha\in J(w)$, the first term is zero and the second is
positive.
For $\alpha\in\Delta\setminus(\Delta_{P_0}^P\cup J(w))$,
the first term is at least one.
A sufficiently small $\varepsilon_0$ therefore gives the second
assertion of \eqref{eq:th45-parameter-properties}.
We may also make $\sigma_T$ dominant.
Let $\overline{\mathcal C}^{\,*}= \left\{
 \lambda\in\mathfrak a^*:
 \langle\lambda,\alpha^\vee\rangle\ge0
 \text{ for every }\alpha\in\Delta
 \right\}$.
Then, $\sC_\rho\cap\overline{\mathcal C}^{\,*}
 = \{ \lambda\in\overline{\mathcal C}^{\,*} \mid
 \rho-\lambda\in \sum_{\alpha\in\Delta}\mathbb R_{\ge0}\alpha\}$.
Equivalently, a dominant element \(\lambda\in\mathfrak a^*\)
belongs to $\sC_\rho$ if and only if $\rho-\lambda$ is a nonnegative real linear combination
of simple roots. 
Since $\rho-\sigma_T=\frac{\varepsilon_0}{1+\|T\|}\chi$ is a nonnegative combination of simple roots, we obtain $\sigma_T\in\sC_\rho$.
\end{proof}

Define
\[
 \Phi_{T,P,w,l}(\lambda)
 \coloneqq \int_{\mathcal D_P(T)}
   (1+\|H_P-T_P\|)^l e^{-\eta_w(\lambda)(H)}\,dH.
\]
The integral converges absolutely and locally uniformly in $\lambda$ on the tube defined by 
\[
\eta_w(\operatorname{Re}\lambda)(\varpi_\alpha^\vee)>0 \quad \text{for all $\alpha\in\Delta\setminus\Delta_{P_0}^P$,}
\]
and hence defines a holomorphic function there.
\begin{lemma}
\label{lem:th45-cone}
Let $\sigma_T$ be as in Lemma~\ref{lem:th45-parameter}. 
For every constant-coefficient differential operator \(D\) of
order \(k\),
\begin{equation}\label{eq:th45-cone-bound}
 \sup_{\xi\in\fa^*}
 |D\Phi_{T,P,w,l}(\sigma_T+i\xi)|
 \ll_{D,l} (1+\|T\|)^{n+l+k}\cJ_{P,G}(T).
\end{equation}
\end{lemma}

\begin{proof}
Put
\[
 \mathcal B_P(T)\coloneqq
 \{Y\in\fa_{P_0}^P\mid
   \tau_{P_0}^P(Y-T_1)=1,\ 
   \widehat\tau_{P_0}^P(T-Y)=1\}.
\]
Up to boundary sets of measure zero, the domain \(\mathcal D_P(T)\)
has the parametrization
\begin{equation}\label{eq:th45-cone-coordinates}
 H=H_b+\sum_{\alpha\in\Delta\setminus\Delta_{P_0}^P}
               t_\alpha\varpi_\alpha^\vee,
 \qquad H_b=Y+T_P,\quad Y\in\mathcal B_P(T),\quad t_\alpha>0.
\end{equation}

For $\alpha\in\Delta_{P_0}^P$, we have $\alpha(H_b)\geq\alpha(T_1)$.
For $\alpha\in\Delta\setminus \Delta_{P_0}^P$, $Y-T^P$ is a negative combination of the coroots in
$\Delta_{P_0}^P$, hence the off-diagonal Cartan inequalities give
\[
 \alpha(H_b)=\alpha(T)+\alpha(Y-T^P)\geq\alpha(T).
\]

Suppose $P^+\neq G$. 
Then, there exists an element $\alpha\in\Delta\setminus \Delta_{P_0}^{P^+}$ such that the $\theta$-orbit of $\alpha$ is disjoint from $\Delta_{P_0}^{P^+}$.
By $\operatorname{sat}_\theta\bigl(\Delta_{P_0}^P\cup I(w)\bigr)=\Delta$, the $\theta$-orbit meets $I(w)$, and hence also
$\theta^{-1}I(w)$.
Equation~\eqref{eq:th45-eta-rho} consequently implies that, for some $c>0$
\[
\eta_w(\rho)(H_b)\geq d(T)-c. 
\]
If $P^+=G$, the same equation gives $\eta_w(\rho)(H_b)\geq-c$. 
Hence, it follows that
\[
 e^{-\eta_w(\rho)(H_b)}\ll\mathcal J_{P,G}(T).
\]
Recall $L_w= 1-w\theta$, so $\eta_w(\lambda)(H)=\lambda(L_wH)$.
Then, \eqref{eq:th45-parameter} and $\|H_b\|\ll (1+\|T\|)$ give
\begin{equation}\label{eq:th45-base-bound}
 e^{-\eta_w(\sigma_T)(H_b)}
 =e^{-\eta_w(\rho)(H_b)}
   e^{\varepsilon_0\chi(L_wH_b)/(1+\|T\|)}
 \ll\mathcal J_{P,G}(T).
\end{equation}

A derivative of order $k$ inserts a polynomial in $L_wH$ bounded
by $O_D((1+\|H\|)^k)$. 
Using \eqref{eq:th45-parameter-properties}, \eqref{eq:th45-cone-coordinates},and \eqref{eq:th45-base-bound}, we obtain
\[
 |D\Phi_{T,P,w,l}(\sigma_T+i\xi)|
 \ll_{D,l}\cJ_{P,G}(T) \, A_T^{n-q}\int_{\R_{>0}^{q}}
 (1+\textstyle\sum_\alpha t_\alpha)^l
 (A_T+\textstyle\sum_\alpha t_\alpha)^k
 e^{-cA_T^{-1}\sum_\alpha t_\alpha}\,dt,
\]
where $A_T\coloneqq 1+\|T\|$ and $q\coloneqq\dim\fa_P^G$.
The substitution \(t_\alpha=A_Tv_\alpha\) proves
\eqref{eq:th45-cone-bound}.
\end{proof}

\begin{lemma}
\label{lem:th45-bruhat}
There is an integer $b\geq0$ such that
\begin{equation}\label{eq:th45-bruhat-bound}
 J_{P,w,l}^T(\tf)
 \ll_l (1+\|T\|)^{n+l+|D_w|}\, \cJ_{P,G}(T)\,  \|\tf\|_{G(\A)^1,b},
 \quad D_w\coloneqq \{\alpha\in\Delta\mid w^{-1}\alpha<0\}.
\end{equation}
\end{lemma}
\begin{proof}
For $\lambda\in\fa^*_\C$,
\[
 \fI(\lambda)
 \coloneqq
 \left\{
 \varphi\colon G(\mathbb A)\longrightarrow\mathbb C
 \text{ smooth}
 \;\middle|\;
 \varphi(pg)=p^{\lambda+\rho}\varphi(g)
 \text{ for all }p\in P_0(\mathbb A),\ g\in G(\mathbb A)
 \right\}.
\]
The intertwining operator $M(w^{-1},\lambda)\colon \fI(\lambda)\to \fI(w^{-1}\lambda)$ is obtained by the meromorphic continuation of the integral
\[
 (M(w^{-1},\lambda)\varphi)(g)
 \coloneqq
 \int_{N_{w^{-1}}(\A)\backslash N_0(\A)}
 \varphi(n_wug)\,du.
\]
This integral converges if $\operatorname{Re}\lambda-\rho$ lies in the positive Weyl chamber. 
Then, the spherical scalar $m_w(\lambda)$ is defined by $(M(w^{-1},\lambda)\varphi)(1)=m_w(\lambda)\varphi(1)$ for a right $K$-invariant function $\varphi\in\fI(\lambda)$.
By \cite{FL11semi}*{Lemma~3.7}, there exists $r_0>0$, depending
only on the fixed data, such that, for every $w\in W$, the function
\[
 h_w(\lambda) = m_w(\lambda)
 \prod_{\alpha\in D_w}
 \langle\lambda-\rho,\alpha^\vee\rangle
\]
is holomorphic and of moderate growth on the tube
\[
 \mathcal T_{r_0}\coloneqq 
 \left\{
 \lambda\in\fa_{\mathbb C}^* \mid
 \|\operatorname{Re}\lambda-\rho\|<r_0
 \right\}.
\]
By Cauchy's estimates, for every constant-coefficient differential
operator $D$, there exists $N_D\ge0$ such that
\begin{equation}\label{eq:Dh_w}
 |D h_w(\sigma+i\xi)|
 \ll_D (1+\|\xi\|)^{N_D},
 \quad
  \|\sigma-\rho\|\le \frac{r_0}{2},
 \quad \xi\in\fa^*,
 \quad w\in W.    
\end{equation}

Let $\Sigma_{\mathrm{red}}^+$ denote the set of positive reduced roots for $P_0$, and put
\[
 \Phi_{w^{-1}} \coloneqq  \{\beta\in\Sigma_{\mathrm{red}}^+\mid w^{-1}\beta<0\}.
\]
Then the rank-one factorization of \(m_w(\lambda)=m(w^{-1},\lambda)\) is
\[
 m_w(\lambda) =
 \prod_{\beta\in\Phi_{w^{-1}}}
 m_\beta\bigl(\langle\lambda,\beta^\vee\rangle\bigr),
\]
where $m_\beta(s)$ is the scalar associated with the spherical intertwining operator for the Levi subgroup $M_\beta$ of relative semisimple rank one corresponding to $\beta$, characterized by $\mathfrak a_{M_0}^{M_\beta}=\mathbb R\beta^\vee$. 
Since every coroot in $\Phi_{w^{-1}}$ lies in the span of the coroots in $I(w)$, 
the annihilation condition in Lemma~\ref{lem:th45-parameter} implies
\[
 m_w(\lambda+z\chi)=m_w(\lambda) \quad (z\in\R),\qquad
 m_w(\lambda)=
 \frac{h_w(\lambda)}
      {\prod_{\alpha\in D_w}
                       \langle\lambda-\sigma_T,\alpha^\vee\rangle}.
\]

Choose a sufficiently small fixed $t_0>0$ and put
$\lambda_0\coloneqq \sigma_T+t_0\rho$.
For $\beta\in\Phi_{w^{-1}}$, 
\[
\langle \lambda_0,\beta^\vee\rangle=\langle (1+t_0)\rho,\beta^\vee\rangle.
\]
Hence, the integral $M(w^{-1},\lambda)$ converges absolutely for $\operatorname{Re}\lambda=\lambda_0$. 
By the Mellin inversion, for some constant $c_0> 0$
\[
F_{\tf}(\exp(H)n_wn)=c_0 \int_{\fa^*} \varphi_{\tf}(n_wn,\lambda_0+i\xi)\, e^{(\lambda_0+i\xi+\rho)(H)}d\xi.
\]
Therefore, by $\varphi_{\tf}(\cdot,\lambda)\in I(\lambda)$ and Fubini's theorem, 
\[
\int_{N_{w^{-1}}(\A)\bs N_0(\A)} F_{\tf}(\exp(H)n_wn)dn=c_0 \, \int_{\fa^*}m_w(\lambda_0+i\xi) \varphi_{\tf}(\lambda_0+i\xi)\, e^{(\lambda_0+i\xi+\rho)(H)}d\xi.
\]
By $2\rho_{w^{-1}}(\theta H)=\rho(L_wH)$, 
\[
e^{2\rho_{w^{-1}}(\theta H)}e^{-(\lambda+\rho)(L_wH)}=e^{-\lambda(L_wH)}=e^{-\eta_w(\lambda)(H)}.
\]
By Lemma \ref{lem:th45-parameter},
\[
 \eta_w(\lambda_0)(\varpi_\alpha^\vee)
 =\eta_w(\sigma_T)(\varpi_\alpha^\vee)
       +t_0\eta_w(\rho)(\varpi_\alpha^\vee)>0
 \quad(\alpha\in\Delta\setminus\Delta_{P_0}^P).
\]
Therefore, we obtain
\begin{equation}\label{eq:th45-mellin}
 J_{P,w,l}^T(\tf)
 =c_0\int_{\operatorname{Re}\lambda=\lambda_0}
       \Phi_{T,P,w,l}(\lambda)m_w(\lambda)
       \varphi_{\tf}(\lambda)\,d\lambda.
\end{equation}

We now apply \cite{FL16}*{Proposition~3.4}, with
\[
 V=\fa,\qquad \mu_0=\sigma_T,\qquad
 S=\{\alpha^\vee\mid \alpha\in D_w\},\qquad
 h(\lambda)=\Phi_{T,P,w,l}(\lambda)h_w(\lambda).
\]
The additive test function is $\Psi_{\tf,\sigma_T}$, for which
\[
 \widehat\Psi_{\tf,\sigma_T}(\lambda-\sigma_T)
 =\varphi_{\tf}(\lambda).
\]
Here the Fourier--Laplace transform is defined using
$e^{-\lambda(H)}$.
The proposition bounds the integral by a fixed continuous seminorm
of this test function and the polynomial-growth bounds for derivatives
of \(h\) on \(\sigma_T+i\fa^*\), through order \(|D_w|\).
Its constant does not involve the width of the holomorphy
neighborhood.
Indeed, after translating \(\sigma_T\) to zero, the pole directions
and the displacement \(\lambda_0-\sigma_T=t_0\rho\) are fixed.
We choose \(\varepsilon_0,t_0\) small enough that the segment between
\(\sigma_T\) and \(\lambda_0\) lies in the fixed tube for \(h_w\);
it also lies in the convergence tube.
By Lemma~\ref{lem:th45-cone} and \eqref{eq:Dh_w}, the proposition consequently gives, for some fixed \(b\),
\[
 J_{P,w,l}^T(\tf)
 \ll_l (1+\|T\|)^{n+l+|D_w|}\cJ_{P,G}(T)
       \sum_{|\mathbf a|\leq b}
       \|\partial_H^{\mathbf a}\Psi_{\tf,\sigma_T}\|_{L^1(\fa)}.
\]
Lemma~\ref{lem:th45-abel} applies because
$\sigma_T\in\sC_\rho$, and proves
\eqref{eq:th45-bruhat-bound}.
\end{proof}

We can now complete the proof of the theorem.
The preceding reduction and \eqref{eq:th45-bruhat-bound} prove the asserted inequality for
any bi-$K$-invariant and non-negative functions in $C^\inf_c(\tG(\A)^1)$.
Let $f\in\cC(\tG(\A)^1)$ and suppose that $f$ is non-negative and bi-$K$-invariant. 
By the same reduction as before and Lemma \ref{lem:1} (3), 
there exists a bi-$K$-invariant function $\tf\in \cC(G(\A)^1)$ such that $\|\tf(x)\|_{G(\A)^1,m} \ll_m  \|f(x)\|_{\tG(\A)^1,m+s_1}$ for some $s_1>0$ and the same bound as \eqref{eq:r3} holds. 
Hence, it suffices to prove $J^T_{P,w,l}(\tf)\ll_l (1+\|T\|)^{r_1} \cJ_{P,G}(T)\|\tf\|_{G(\A)^1,s_2}$ for some $r_1\ge 0$ and some $s_2>0$. 
There is a sequence $\{\tf_j\}_{j=1}^\inf$ in $C^\inf_c(G(\A)^1)$ such that $\tf_j$ converges to $\tf$ in the topology of $\cC(G(\A)^1)$, and $\tf_j$ is non-negative and bi-$K$-invariant. 
By Lemma \ref{lem:th45-bruhat} and Fatou's lemma, for some $r_2\ge0$ and some $s_3>0$, 
\begin{align*}
J^T_{P,w,l}(\tf)\le&  \liminf_{j\to \inf} J^T_{P,w,l}(\tf_j)\ll_l (1+\|T\|)^{r_2} \cJ_{P,G}(T) \liminf_{j\to \inf} \|\tf_j\|_{G(\A)^1,s_3} \\
\le & (1+\|T\|)^{r_2} \cJ_{P,G}(T) \|\tf\|_{G(\A)^1,s_3}.    
\end{align*}
This proves both the absolute convergence and the required bound.
\end{proof}

\begin{remark}\label{rem:compareFL}
In the proof of \cite{FL16}*{Proposition~3.2}, Finis and Lapid use two real
parameters. The first, $\lambda_0$, specifies the initial contour of Mellin
inversion and is chosen so that $\lambda_0-\rho_0$ lies in the open positive
Weyl chamber. The second, $\mu_0$, specifies the real part at which the
derivatives are estimated in \cite{FL16}*{Proposition~3.4}.
In the untwisted setting, the condition $Q(w)P=G$ implies that
$\rho_0-w^{-1}\rho_0$ is positive on every unbounded ray of the integration
domain, as shown in \cite{FL11semi}*{Lemma 2.2}. 
Thus one may take $\mu_0=\rho_0$ and choose both parameters independently of $T$.

In the twisted setting, the corresponding exponent is
$\eta_w(\rho)=\theta^{-1}(\rho-w^{-1}\rho)$, whose simple-root support is
$\theta^{-1}I(w)$. The condition
$\operatorname{sat}_{\theta}(\Delta_{P_0}^{P}\cup I(w))=\Delta$
allows $\eta_w(\rho)$ to vanish on some unbounded rays of $D_P(T)$.
For example, in split type $A_2$, if $\theta$ interchanges the two simple
roots, $P=P_0$, and $w=s_{\alpha_1}$, then $\eta_w(\rho)=\alpha_2$,
which vanishes on the ray generated by $\varpi_{\alpha_1}^{\vee}$.
Replacing $\rho$ by a positive multiple of $\rho$ leaves this vanishing
unchanged.

We therefore take $\mu_0=\sigma_T=\rho-\varepsilon_0\chi/A_T$, where
$A_T=1+\|T\|$. The size of this perturbation serves two purposes.
On each unbounded ray, we have
$\eta_w(\sigma_T)(\varpi_\alpha^\vee)\geq c/A_T$, so integration in
the ray parameter gives at most a polynomial loss in $A_T$. On the
bounded base of $\cD_P(T)$, we have $\|H_b\|\ll A_T$, and hence $e^{-\eta_w(\sigma_T)(H_b)}
=e^{-\eta_w(\rho)(H_b)}
 \exp\!\left(\frac{\varepsilon_0}{A_T}\chi((1-w\theta)H_b)\right)
\ll e^{-\eta_w(\rho)(H_b)}$. 
Thus the perturbation preserves the bound $\cJ_{P,G}(T)$ on the base.
Moreover, $\sigma_T\in \sC_\rho$, so the Abel-transform estimates are
uniform in $T$.

Finally, we take $\lambda_0=\sigma_T+t_0\rho$ with $t_0>0$ fixed.
Since $\chi$ annihilates $\beta^\vee$ for every positive reduced root
$\beta$ such that $w^{-1}\beta<0$, we have
$\langle\lambda_0,\beta^\vee\rangle
=(1+t_0)\langle\rho,\beta^\vee\rangle$ for every such $\beta$.
Thus the initial contour follows the evaluation point $\sigma_T$ at
a fixed displacement, while the real parameters governing absolute
convergence of the intertwining integral remain independent of $T$.
\end{remark}

\subsection{Proof of Theorem \ref{thm:main theorem 1} (2)}
Set $\sU\coloneqq\{T\in \fa_{M_0}^G \mid d(T)>d_*\}$. 
For $T\in\mathscr U$, let
\[
 \mathbf J_T(f)
 \coloneqq \bigl(J_{\tfo}^{T,\tG}(f)\bigr)_{\tfo\in\tcO}
 \in\ell^1(\tcO).
\]
By Theorem \ref{thm:main theorem 1} (1), proved above,
\begin{equation}\label{eq:gr-polynomial-evaluation-bound}
 \|\mathbf J_T(f)\|_1
 \le C\level(K_0)^s(1+\|T\|)^r
       \|f\|_{\tG(\A)^1,m}.
\end{equation}
In particular, $\mathbf J_T$ is a continuous linear map from
$\cC(\tG(\A)^1;K_0)$ to $\ell^1(\tcO)$.

Suppose first that $f\in C_c^\infty(\tG(\A)^1;K_0)$.
We use the geometric translation argument in the proof of
\cite{LW13}*{Theorem~11.1.1}.
For $\tQ\in\cF^{\tG}(\tM_0)$ with $\tQ\supset\tP_0$, put
\[
 \gamma_{\tQ}(X)
 \coloneqq \int_{\fa_{\tQ}^{\tG}}\Gamma_{\tQ}(H,X)\,dH.
\]
These functions $\Gamma_{\tQ}$ and $\gamma_{\tQ}$ are defined in \cite{LW13} \S~2.9, p.47, and \S~11.1, p.~149.
By \cite{LW13}*{Lemmas~1.9.1 and~2.9.2}, this is a homogeneous
polynomial of degree $\dim\fa_{\tQ}^{\tG}$, viewed as a polynomial
on $\fa_{M_0}^G$ through the natural projection. In particular,
\[
 \gamma_{\tG}\equiv1,\qquad
 \gamma_{\tQ}(0)=0\quad(\tQ\ne\tG).
\]

Fix $T\in\mathscr U$. The cone identity
\cite{LW13}*{Proposition~2.9.4~(2)} gives, for all sufficiently positive $X\in \fa_{M_0}^G$,
\begin{equation}\label{eq:gr-polynomial-translation}
 J_{\tfo}^{T+X,\tG}(f)
 =\sum_{\substack{\tQ\in\cF^{\tG}(\tM_0)\\Q\supset P_0}}
   \gamma_{\tQ}(X)A_{\tQ,\tfo}^{T}(f).
\end{equation}
Here $A_{\tQ,\tfo}^{T}(f)$ denotes the descended Levi integral,
with all Levi coarse classes contributing to $\tfo$ grouped
together. It is independent of $X$, and $A_{\tG,\tfo}^{T}(f)=J_{\tfo}^{T,\tG}(f)$.

For each $T\in\sU$, define the polynomial
\[
 P_{\tfo,f;T}(\bT)
 \coloneqq \sum_{\substack{\tQ\in\cF^{\tG}(\tM_0)\\Q\supset P_0}}
   \gamma_{\tQ}(\bT-T)A_{\tQ,\tfo}^{T}(f)
 \qquad(\bT\in \fa_{M_0}^G).
\]
Its degree is at most $\dim\fa_{M_0}^G$.
For any $T,T'\in\sU$, \eqref{eq:gr-polynomial-translation} shows that
$P_{\tfo,f;T}(\bT)$ and $P_{\tfo,f;T'}(\bT)$ both equal
$J_{\tfo}^{\bT,\tG}(f)$ whenever $\bT-T$ and $\bT-T'$ are sufficiently
positive. 
They therefore agree on a nonempty open subset of $\fa_{M_0}^G$
and hence are identical.
Denote their common value by $P_{\tfo,f}$.
Evaluating the defining polynomial at $\bT=T$ gives
\begin{equation}\label{eq:gr-polynomial-compact}
 J_{\tfo}^{T,\tG}(f)=P_{\tfo,f}(T)
 \qquad(T\in\sU).
\end{equation}
Thus the equality holds on the fixed region $\sU$,
independently of the support of $f$.

Let $\cP_{\le n}(\fa_{M_0}^G)$ denote the set of polynomials on $\fa_{M_0}^G$, whose degrees are less than $n+1$. 
Fix a monomial basis $(p_\alpha)_{\alpha\in\mathscr A}$ of
$\cP_{\le n}(\fa_{M_0}^G)$ and put $N\coloneqq |\mathscr A|$.
Since $\sU$ is nonempty and open, choose
$T^{(1)},\ldots,T^{(N)}\in\sU$ such that evaluation
at these points is an isomorphism
$\mathcal P_{\le n}(\fa_{M_0}^G)\longrightarrow\C^N$.
Such points exist because a polynomial vanishing on
$\mathscr U$ is zero.
Let $b_{\alpha i}$ be the interpolation coefficients, so that
\[
 \mathbf{P}(T)=\sum_{\alpha\in\mathscr A}p_\alpha(T)
       \sum_{i=1}^{N}b_{\alpha i}\mathbf{P}(T^{(i)})
 \qquad(\mathbf{P}\in\mathcal P_{\le n}(\fa_{M_0}^G)).
\]

For every $f\in\cC(\tG(\A)^1;K_0)$, define
\[
 \mathbf c_\alpha(f)
 \coloneqq\sum_{i=1}^{N}b_{\alpha i}\mathbf J_{T^{(i)}}(f)
 =\bigl(c_{\tfo,\alpha}(f)\bigr)_{\tfo\in\tcO}.
\]
These are continuous linear maps with values in $\ell^1(\tcO)$.
Moreover, \eqref{eq:gr-polynomial-evaluation-bound} gives
\begin{equation}\label{eq:gr-polynomial-coefficients}
 \sum_{\tfo\in\tcO}|c_{\tfo,\alpha}(f)|
 =\|\mathbf c_\alpha(f)\|_1  \ll   \level(K_0)^s\|f\|_{\tG(\A)^1,m}.
\end{equation}

For compactly supported $f$, interpolation and
\eqref{eq:gr-polynomial-compact} imply
\[
 \mathbf J_T(f)
 =\sum_{\alpha\in\mathscr A}p_\alpha(T)\mathbf c_\alpha(f)
 \qquad(T\in\mathscr U).
\]
For each fixed $T\in\sU$, both sides are continuous
linear maps from $\cC(\tG(\A)^1;K_0)$ to $\ell^1(\tcO)$.
Hence, the density of $C_c^\infty(\tG(\A)^1;K_0)$ extends
this identity to every $f\in\cC(\tG(\A)^1;K_0)$.
Consequently,
\begin{equation}\label{eq:gr-polynomial-continuation}
 J_{\tfo}^{T,\tG}(f)
 =\sum_{\alpha\in\mathscr A}c_{\tfo,\alpha}(f)p_\alpha(T)
 \qquad(T\in\mathscr U),
\end{equation}
and the right-hand side defines the unique polynomial
continuation to all $T\in \fa_{M_0}^G$, of degree at most $n$.
Its coefficients are continuous linear forms and satisfy
\eqref{eq:gr-polynomial-coefficients}.
By summing the coefficients as in $c_\alpha(f)=\sum_{\tfo\in\tcO}c_{\tfo,\alpha}(f)$ and applying Theorem~\ref{thm:main theorem 1} (1), we obtain the analogous result for $J^{T,\tG}(f)$.
This completes the proof of Theorem~\ref{thm:main theorem 1}~(2).

\section{On the continuity of the geometric expansion for refined classes}
\label{fg:section}

We retain the notation and measures of \S~\ref{sec:preliminaries}, \S~\ref{sec:cga}, and \S~\ref{gr:section}.

\subsection{Statement of the theorem}\label{sec:mainrefinedclass}

For $\gamma$, $\gamma'\in\tG(\Q)$, write $\gamma\sim\gamma'$
if their semisimple parts are conjugate under $G(\Q)$ and
$\gamma$ and $\gamma'$ are conjugate under $G(\overline\Q)$.
Let $\cO$ denote the set of these refined classes, and let
$\tfo(\fo)\in\tcO$ be the coarse class containing $\fo\in\cO$.

We recall the notion of induction of conjugacy classes, see \cites{LS79,Hof16}. Our presentation follows \cite{FL16}*{\S~6.3}.
For $\tP=\tM_PN_P\in\cF^{\tG}(\tM_0)$ with $P\supset P_0$ and
$\gamma\in\tM_P(\Q)$, denote by
$\cI_{\tP}^{\tG}(\gamma)$ the unique geometric
$G$-conjugacy class whose intersection with $\gamma N_P$ is open and dense.
Define
\begin{equation}\label{fg:induction}
 I_{\tP}^{\tG}(\gamma)
 \coloneqq\mathcal I_{\tP}^{\tG}(\gamma)(\Q)\cap\tfo_\gamma\in\cO,
\end{equation}
where $\tfo_\gamma$ is the coarse class containing $\gamma$. 
See Lemma \ref{lem:well defined} for the well definedness of \eqref{fg:induction}. 
For general $\gamma\in\tP$, we use the same notation for the induced class of its Levi projection to $\tM_P$.

For $f\in\cC(\tG(\A)^1;K_0)$ and $\fo\in\cO$, define
\begin{align}
 K_{\tP,\fo}(x,x)
 &\coloneqq\sum_{\substack{\gamma\in\tM_P(\Q)\\
                    I_{\tP}^{\tG}(\gamma)=\fo}}
       \int_{N_P(\A)}f(x^{-1}\gamma n x)\,dn,
       \label{fg:parabolic-kernel}\\
 k_\fo^T(x,f)
 &\coloneqq\sum_{\substack{\tP\in\cF^{\tG}(\tM_0)\\P\supset P_0}}
       (-1)^{\dim\fa_{\tP}^{\tG}}
       \sum_{\xi\in P(\Q)\backslash G(\Q)}
       K_{\tP,\fo}(\xi x,\xi x;f)
       \widehat\tau_{\tP}(H_P(\xi x)-T).
       \label{fg:truncated-kernel}
\end{align}
The positive sum obtained by taking absolute values in
\eqref{fg:parabolic-kernel} is bounded by the unrestricted
parabolic sum. Lemma~\ref{lem:sup kernel} therefore gives
absolute convergence at every $x$, as well as the local
finiteness needed in \eqref{fg:truncated-kernel}.

Define $B_{\fo,Q,R}^T$ by replacing the coarse parabolic
kernels in $B_{\tfo,Q,R}^T$ by
\eqref{fg:parabolic-kernel}, and put
\[
 m_\fo^T(x,f)
 \coloneqq\sum_{\substack{Q\subset R\\Q^+=R^-}}
       B_{\fo,Q,R}^T(x,f).
\]
Thus every one-term interval is retained in the main kernel.
\eqref{eq:lw911-sigma-identity} and Lemma \ref{lem:partitionformula} give
\[
 k_\fo^T(x,f)=\sum_{Q\subset R}B_{\fo,Q,R}^T(x,f)
 \qquad(d(T)>d_0).
\]

\begin{theorem}\label{thm:main theorem 2}
There exist constants $C,d_*>0$ and integers $m,s,r\ge0$,
with $d_*>d_0$, such that the following assertions hold for
every compact open subgroup $K_0\subset G(\A_f)$ and every
$f\in\cC(\tG(\A)^1;K_0)$.
\begin{enumerate}[(1)]
\item For every $T\in\fa_{M_0}^G$ satisfying $d(T)>d_*$,
\begin{equation*}
 \sum_{\fo\in\cO}\int_{G(\Q)\bs G(\A)^1}
     |k_\fo^T(x,f)|dx
 \le C\level(K_0)^s
          (1+\|T\|)^r\|f\|_{\tG(\A)^1,m}.
\end{equation*}
In particular, the integrals $J_\fo^{T,\tG}(f)$ $(\fo\in\cO)$ are absolutely convergent, and we have
\[
 J^{T,\tG}(f)=\sum_{\fo\in\cO}J_\fo^{T,\tG}(f).
\]
\item $J_\fo^{T,\tG}(f)$ $(\fo\in\cO)$ are polynomials
in $T$ of degree at most $\dim\fa_{M_0}^G$.
Their coefficients are continuous linear forms in $f$.
For a fixed monomial basis, the sum over $\fo$ of the absolute
values of each coefficient is bounded by
$C\level(K_0)^s\|f\|_{\tG(\A)^1,m}$,
after increasing $C$ if necessary.
\item For every $T\in\fa_{M_0}^G$ satisfying $d(T)>d_*$, the approximation estimate
\begin{align*}
 \sum_{\fo\in\cO}|J_\fo^{T,\tG}(f)-M_\fo^T(f)|
 &\le\sum_{\fo\in\cO}\int_{G(\Q)\bs G(\A)^1}
           |k_\fo^T(x,f)-m_\fo^T(x,f)|\,dx\\
 &\le C\level(K_0)^s
       (1+\|T\|)^r e^{-d(T)}\|f\|_{\tG(\A)^1,m}
\end{align*}
holds. Note that $\sum_{\fo\in\cO}\int_{G(\Q)\bs G(\A)^1}|m_\fo^T(x,f)|\,dx\le C\level(K_0)^s
       (1+\|T\|)^r\|f\|_{\tG(\A)^1,m}$.
\end{enumerate}
\end{theorem}

The refined and coarse kernels satisfy
\begin{equation}\label{fg:coarse-compatibility}
 \sum_{\fo:\tfo(\fo)=\tfo}k_\fo^T=k_{\tfo}^T,
 \qquad
 \sum_{\fo:\tfo(\fo)=\tfo}m_\fo^T=m_{\tfo}^T.
\end{equation}
Indeed, each rational Levi element has exactly one induced
refined class, contained in its coarse class.
The sums in \eqref{fg:coarse-compatibility} are finite by
\eqref{fg:finite-refinement} below. Consequently, 
\[
 J_{\tfo}^{T,\tG}(f)
 =\sum_{\fo:\tfo(\fo)=\tfo}J_\fo^{T,\tG}(f),
 \qquad
 M_{\tfo}^T(f)
 =\sum_{\fo:\tfo(\fo)=\tfo}M_\fo^T(f).
\]
In the following subsections, we prove Theorem \ref{thm:main theorem 2}.

\subsection{Algebraic properties of induction}

There is a constant $b_G$ such that
\begin{equation}\label{fg:finite-refinement}
 \#\{\fo\in\cO:\tfo(\fo)=\tfo\}\le b_G
 \qquad(\tfo\in\tcO).
\end{equation}
Indeed, for a fixed semisimple representative $\sigma\in\tG(\Q)$, the
classes in question are indexed by the orbits of the
centralizer $Z_G(\sigma)$ of $\sigma$ in $G$ on its unipotent elements that
contain rational points.
Every such unipotent element belongs to $Z_G(\sigma)^\circ$, which is the identity component of $Z_G(\sigma)$.
The number of its geometric unipotent classes is bounded
uniformly in $\sigma$, since its rank and dimension are bounded
in terms of $G$.
Hence, such a constant $b_G$ exists. 

\begin{lemma}\label{lem:well defined}
Let $\tP=\tM N_P$ and $\gamma\in\tM(\Q)$. 
The geometric conjugacy class $I_{\tP}^{\tG}(\gamma)$ is defined over $\Q$ and has rational points. 
Further, $(I_{\tP}^{\tG}(\gamma)\cap \gamma N_P)(\Q)$ is contained in $\tfo_\gamma$. 
In particular, \eqref{fg:induction} is well defined.  
\end{lemma}
\begin{proof}
The conjugation argument of \cite{FL16}*{Lemma~6.7} applies to
$\gamma N_P$: on the vector quotients of a suitable invariant
central filtration of $N_P$, separate the fixed subspace of
$\gamma_s$ from its complement.
On that complement $\Ad(\gamma)-1$ is invertible, since
$\gamma_s$ acts semisimply and $\gamma_u$ acts unipotently
and commutes with $\gamma_s$.
Solving successively on these quotients conjugates an element
of $\gamma N_P$ into $\gamma N_{P,\gamma_s}$.
The construction is defined over $\Q$ when $\gamma$ is rational.
It follows that the semisimple part of every element of
$\gamma N_P$ is $N_P$-conjugate to $\gamma_s$, with a rational
conjugating element for rational points.
There are consequently only finitely many geometric $G$-orbits
meeting $\gamma N_P$, because there are finitely many
geometric unipotent classes in $Z_G(\gamma_s)^\circ$.
The irreducibility of $\gamma N_P$ and local closedness of
orbits give the unique open dense intersection.
It is defined over $\Q$ and has rational points, since
$N_P$ is a rational affine space over $\Q$.
Thus the proof is completed.        
\end{proof}

If $\tP\subset\tR$
and $\eta\in\mathcal I_{\tP\cap\tM_R}^{\tM_R}(\gamma)$,
then $\mathcal I_{\tR}^{\tG}(\eta) =\mathcal I_{\tP}^{\tG}(\gamma)$.
This follows by comparing the dense orbits in the closures
obtained from $N_P=(N_P\cap M_R)N_R$.
Likewise, if $\gamma'\in\overline{M\cdot\gamma}$, then
\[
 \mathcal I_{\tP}^{\tG}(\gamma')
   \subseteq\overline{\mathcal I_{\tP}^{\tG}(\gamma)}.
\]
In particular, for $\gamma\in\tP$ and $\tP\subset\tR$,
\[
 \mathcal I_{\tR}^{\tG}(\gamma)
   \subseteq\overline{\mathcal I_{\tP}^{\tG}(\gamma)}.
\]
These are the properties corresponding to
\cite{FL16}*{(27)--(29)}.

\begin{lemma}\label{fg:exceptional-polynomial}
There is an integer $d_G$ such that, for every standard
$\tP=\tM N_P$ and every $\gamma\in\tM(\Q)$, there is a
nonzero polynomial $F_\gamma\in\Q[\Lie(N_P)]$ of degree
at most $d_G$ satisfying
\[
 \gamma\exp X\notin\mathcal I_{\tP}^{\tG}(\gamma)
 \quad\Longrightarrow\quad F_\gamma(X)=0.
\]
The same degree bound can be used for induction inside any
Levi subset.
\end{lemma}
\begin{proof}
We give a direct proof of the twisted analogue of
\cite{FL16}*{Lemma~6.6}, using the dimensions of conjugacy
orbits and minors of the adjoint action.

Put $\fO\coloneqq \cI_{\tP}^{\tG}(\gamma)$, $r\coloneqq \dim\fO$, and $d\coloneqq \dim G$. 

Since $\fO\cap\gamma N_P$ is Zariski dense in
$\gamma N_P$, we have
\[
 \gamma N_P\subset\overline{\fO}.
\]
The orbit $\fO$ is locally closed, hence open in its closure.
Since $G$ is connected, $\fO$ and its closure are irreducible.
Therefore, the boundary $\overline{\fO}\setminus\fO$
is a proper closed subset of $\overline{\fO}$ and has
dimension strictly smaller than $\dim\fO$.

In characteristic zero, the stabilizer of $x$ under
conjugation is smooth and has Lie algebra
$\ker(\Ad(x)-1\mid\fg)$.
Thus
\[
 \dim(G\cdot x)=\rank(\Ad(x)-1\mid\fg).
\]
It follows that, for $x\in\gamma N_P$,
\[
 x\notin\fO
 \quad\Longleftrightarrow\quad
 \rank(\Ad(x)-1\mid\fg)<r.
\]
If $r=0$, the exceptional set is empty and we take
$F_\gamma=1$.

Suppose that $r>0$.
Choose $\Q$-bases of $\fg$ and $\Lie(N_P)$.
For $X\in\Lie(N_P)$, the endomorphism $\ad X$ of $\fg$
is nilpotent, so $(\ad X)^d=0$ and
\[
 A_\gamma(X)
 \coloneqq\Ad(\gamma\exp X)-1
 =\Ad(\gamma)\sum_{j=0}^{d-1}
       \frac{(\ad X)^j}{j!}-1.
\]
Since $\gamma$ is rational, the entries of $A_\gamma(X)$
are polynomials over $\Q$ of degree at most $d-1$.
Its maximal rank is $r$, attained whenever
$\gamma\exp X\in\fO$.
We may therefore choose an $r\times r$ minor
$F_\gamma(X)$ that is not the zero polynomial.
The rank characterization above gives
\[
 \gamma\exp X\notin\fO
 \quad\Longrightarrow\quad F_\gamma(X)=0,
\]
\[
 \deg F_\gamma\le r(d-1)\le d(d-1).
\]

For induction inside a Levi subset $\tL$, apply the same
argument to the conjugation action of $L$ on $\tL$,
using $\Lie(L)$ in place of $\fg$.
Since $\dim L\le d$, the resulting degree bound is also
at most $d(d-1)$.
Thus we may take $d_G=d(d-1)$ uniformly.
\end{proof}

\subsection{Restricted alternating sums}

Let $Q\subset R$ be standard parabolic subgroups with
$Q^+\subset R^-$, and let $\tS$ satisfy
$\tQ^+\subset\tS\subset\tR^-$.
With $\Gamma_{Q,R;S}$ as in \eqref{eq:tfl47-gamma}, define
\begin{equation}\label{fg:normalized-gamma}
 \Gamma^{\mathrm{ref}}_{Q,R;S}(T,H)
   \coloneqq \begin{cases}
      \frac{1}{\dim\fa_{\tS}^{\tR^-}}\Gamma_{Q,R;S}(T,H)&\text{if $\tS\subsetneq \tR^-$},\\
      0& \text{if $\tS= \tR^-$}.
     \end{cases}
\end{equation}

Let $\widetilde O\subset\tM_S(\Q)$ be invariant under
conjugation by $Q(\Q)\cap M_S(\Q)$.
For $\fo\in\cO$ and $f\in\cC(\tG(\A)^1;K_0)$, define
\begin{equation}\label{fg:restricted-sum}
 \begin{split}
 \mathcal E^{\mathrm{ref}}_{\fo,Q,R;S}
       (g,\widetilde O,f)
  \coloneqq {}&\sum_{\tS\subset\tP\subset\tR^-}
        (-1)^{\dim\fa_{\tP}^{\tG}}
       \sum_{\eta\in\widetilde O}
       \sum_{\substack{\nu\in N_S^{M_P}(\Q)\\
               I_{\tP}^{\tG}(\eta\nu)=\fo}}
       \int_{N_P(\A)}f(g^{-1}\eta\nu n g)\,dn .
 \end{split}
\end{equation}

\begin{lemma}\label{fg:restricted-estimate}
There are an integer $s\ge0$ and finitely many
$X_1,\ldots,X_m\in\cU(\fg_\C^1)$ such that
\begin{equation*}
 \sum_{\fo\in\cO}
   |\mathcal E^{\mathrm{ref}}_{\fo,Q,R;S}
          (g,\widetilde O,f)|\ll \level(K_0)^s
 e^{-\Gamma^{\mathrm{ref}}_{Q,R;S}(T,H_0(g))}
 \sum_{\eta\in\widetilde O}\sum_{j=1}^m
       \int_{N_S(\A)}|(f*X_j)(g^{-1}\eta n g)|\,dn
\end{equation*}
whenever $d(T)>d_0$ and
$F^Q(g,T)\tau_Q^{R^-}(H_Q(g)-T)=1$.
The differential operators, the exponent
$s$, and the implied constant are independent of $\tO$, $K_0$, $f$, $T$, and $g$.
\end{lemma}
\begin{proof}
We follow \cite{FL16}*{Proposition~6.9 and
Corollaries~6.10--6.11}, recording only the modifications.
The transitivity and closure properties established above
justify its Levi reduction and decomposition into monotone
families.

For a root coordinate of weight $\beta_i$, put
\[
 I_{\tomega}=\{i\mid \operatorname{supp}(\beta_i)\cap \tomega\ne\emptyset\},
 \qquad \tomega=\{ \theta^j(\omega)\mid j\in\Z \}\in\Delta_{\tP_0}^{\tR^-}\setminus\Delta_{\tP_0}^{\tS}.
\]
These coordinate sets may overlap, as permitted in
\cite{FL16}*{Lemma~6.1 and Corollary~6.2}.
Fix $\tomega\in\Delta_{\tP_0}^{\tR^-}\setminus\Delta_{\tP_0}^{\tS}$.
If no minimal member of the monotone family contains $\tomega$,
their cancellation argument gives the gain
$e^{-\min_{\alpha\in \tomega}\alpha(H_a)}$, where
$a=\exp H_a$ is the reduced torus variable.
Otherwise, delete the entire orbit $\tomega$ from a minimal member containing it.
As in \cite{FL16}*{(34)}, minimality and transitivity force
the corresponding rational variable into an exceptional
induction locus.
Lemma~\ref{fg:exceptional-polynomial} places this locus in
the zero set of a nonzero polynomial of uniformly bounded
degree.
Every root of the resulting unipotent group has support
meeting $\tomega$; the fixed lower bounds on the simple-root
coordinates therefore give a contraction factor at most
$C e^{-\min_{\alpha\in \tomega}\alpha(H_a)}$.
Applying \cite{FL16}*{Lemmas~6.4 and~4.4} gives the same gain.
This replaces the use of Corollary~6.5 there, since deleting
$\tomega$ need not give an ordinary maximal parabolic inclusion.
The auxiliary unipotent sums are handled by
\cite{FL16}*{Lemma~3.3}, exactly as in that proof.

Fix $\fo\in\cO$.
Since induction preserves the rational semisimple class,
we may replace $\widetilde O$ by
$\widetilde O\cap\tfo(\fo)$.
Subtract the closed-set restricted sums for the closure
of the geometric orbit underlying $\fo$ and for its boundary,
and apply the triangle inequality, as in
\cite{FL16}*{Corollary~6.11}.
This imposes the condition
$I_{\tP}^{\tG}(\eta\nu)=\fo$.

The restricted set of $\eta$ retains the required conjugation
invariance.
Hence, after summing over $\eta$, the reduction in
Proposition~\ref{prop:tfl47} applies.
Its root comparison and compact-conjugation estimates,
together with Lemma~\ref{lem:levelbound}, replace the gain
for $\tomega$ by
\[
 \exp\!\left(
   -\min_{\alpha\in \tomega}
       \{\alpha(T)+\alpha_Q(H_0(g)-T)\}
 \right),
\]
with a factor $\level(K_0)^s$ and finitely many right
differential operators.
These choices are uniform in $\fo$ and all the parameters
of the lemma, since only the bounded polynomial degree
enters the exceptional-set estimate.

Choose a common differential family for all
$\tomega\in\Delta_{\tP_0}^{\tR^-}\setminus\Delta_{\tP_0}^{\tS}$.
If $\tS\subsetneq \tR^-$, by taking the geometric mean of these bounds for $\tomega\in\Delta_{\tP_0}^{\tR^-}\setminus\Delta_{\tP_0}^{\tS}$, we obtain
$e^{-\Gamma^{\mathrm{ref}}_{Q,R;S}(T,H_0(g))}$.
The positive majorant for a fixed $\fo$ only involves
$\eta\in\widetilde O\cap\tfo(\fo)$.
Thus summation over $\fo$ costs at most the factor $b_G$
in \eqref{fg:finite-refinement}, which is absorbed into
the implied constant.
If $\tS=\tR^-$, there is a single parabolic term, and the
assertion follows directly from the triangle inequality,
including the identity operator in the differential family.

All rearrangements are justified by Lemma~\ref{lem:sup kernel}.
\end{proof}

\begin{lemma}\label{fg:cone-integral}
For $y\in M_S(\A)^1$, $H^-\in(\fa_S^G)^{1-\theta}$,
and $d(T)>d_0$, we have
\begin{align*}
 &\int_{\fa_{\tS}^{\tG}}
   \widetilde\sigma_Q^R(H_0(y)+H^-+H^\theta-T)
   e^{-\Gamma^{\mathrm{ref}}_{Q,R;S}(T,H_0(y)+H^-+H^\theta)}\,dH^\theta \\
 &\ll
   \cJ_{S,R^-}(T)\tau_Q^S(H_0(y)-T)
   \bigl(1+\|T\|+\|H_Q^S(y)-T_Q^S\|+\|H^-\|\bigr)^{\dim\fa^{\tG}_{\tR^-}}.
\end{align*}
\end{lemma}
\begin{proof}
If $S=R^{-}$, then $\Gamma^{\mathrm{ref}}_{Q,R;S}=0$, and the assertion follows from the volume bound provided by Lemma \ref{lem:cga-cone-fibers}. We may therefore assume that $S\subsetneq R^{-}$.
If we use the coordinates $u_{\tbeta}$, $\tbeta\in\Delta_{\tP_0}^{\tR^-}\setminus\Delta_{\tP_0}^{\tS}$ from the proof of
Lemma~\ref{lem:analogueFL48}, then
\begin{equation*}
 e^{-\Gamma^{\mathrm{ref}}_{Q,R;S}(T,H_0(y)+H^-+H^\theta)}
 \le e^{-d(T)}e^{-\frac{1}{d^R_S}\sum_{\tbeta} u_{\tbeta}} ,\quad d^R_S=\dim\fa_{\tS}^{\tR^-}.
\end{equation*}
Therefore, the assertion can be proved by the same argument as in the proof of Lemma~\ref{lem:analogueFL48}. 
\end{proof}

\subsection{Proof of Theorem~\ref{thm:main theorem 2}}

Let $f\in\cC(\tG(\A)^1;K_0)$. 
Define $k_{\fo,Q,R}(x,f)$ and $\cR^{\mathrm{ref}}_{Q,R}(T,f)$
as in Section~\ref{sec:reduction(1)(3)}, replacing coarse
classes and kernels by their refined counterparts.
The strata used there give
\begin{equation*}\label{fg:strata}
 k_{\fo,Q,R}(x,f)
 =\sum_{\tQ^+\subset\tS\subset\tR^-}
       \mathcal E^{\mathrm{ref}}_{\fo,Q,R;S}
          (x,\tM_S(\Q)_Q^\circ,f).
\end{equation*}
Their $Q(\Q)\cap M_S(\Q)$-invariance allows us to apply
Lemma~\ref{fg:restricted-estimate}.
We now use the reduction of Section~\ref{sec:reduction(1)(3)},
with Lemma~\ref{fg:restricted-estimate} in place of Proposition~\ref{prop:tfl47}
and Lemma~\ref{fg:cone-integral} in place of
Lemma~\ref{lem:analogueFL48}.
Thus, Theorem~\ref{thm:main theorem 2} (1) and (3) follow from Theorem~\ref{thm:bound 1}. 

Theorem~\ref{thm:main theorem 2} (2) can be proved by the same argument as in the proof of Theorem \ref{thm:main theorem 1} (2) together with Theorem~\ref{thm:main theorem 2} (1) and (3). 
We therefore omit the details.

\bibliography{bibliography}
\bibliographystyle{amsplain}

\end{document}